\documentclass{article}
\usepackage{amsmath,  amsfonts, amsthm, latexsym, amssymb}
\usepackage{a4wide} 
\usepackage{tikz}\usetikzlibrary{patterns}

 \numberwithin{equation}{section}

\newtheorem{theorem}{Theorem}[section]
\newtheorem{definition}[theorem]{Definition}
\newtheorem{proposition}[theorem]{Proposition}
\newtheorem{remark}[theorem]{Remark}
\newtheorem{lemma}[theorem]{Lemma}
 
\newtheorem{corollary}[theorem]{Corollary}

\usepackage[sc]{mathpazo}

\usepackage{color}

\newcommand{\cU}{{\cal U}}

 \newcommand{\Z}{{\mathbb Z}}
  \newcommand{\T}{{\mathbb T}}

\newcommand{\R}{{\mathbb R}}

\newtheorem*{theorem*}{Theorem}

\def\Xint#1{\mathchoice
{\XXint\displaystyle\textstyle{#1}}%
{\XXint\textstyle\scriptstyle{#1}}%
{\XXint\scriptstyle\scriptscriptstyle{#1}}%
{\XXint\scriptscriptstyle\scriptscriptstyle{#1}}%
\!\int}
\def\XXint#1#2#3{{\setbox0=\hbox{$#1{#2#3}{\int}$ }
\vcenter{\hbox{$#2#3$ }}\kern-.6\wd0}}

\def\dashint{\Xint-}

\titlepage
\title{Stationary equilibria and their stability in a Kuramoto MFG with strong interaction}
\author{Annalisa Cesaroni and Marco Cirant}
\date{ }

\begin{document}

\maketitle
\begin{abstract}   
Recently, R. Carmona, Q. Cormier, and M. Soner proposed a Mean Field Game (MFG) version of the classical Kuramoto model, which describes synchronization phenomena in a large population of ``rational'' interacting oscillators. The MFG model exhibits several stationary equilibria, but the characterization of these equilibria and their ability to capture dynamic equilibria in long time remains largely open.

In this paper, we demonstrate that, up to a phase translation, there are only two possible stationary equilibria: the incoherent equilibrium and the self-organizing equilibrium, given that the interaction parameter is sufficiently large. Furthermore, we present some local stability properties of the self-organizing equilibrium.
\medskip

\noindent
{\footnotesize \textbf{AMS-Subject Classification}}. {\footnotesize  
35Q89, 
49N80, 
92B25 
 }
\end{abstract} 
 

\section{Introduction} 

The classical Kuramoto model is a system of nonlinear ordinary differential equations that describes the dynamics of coupled oscillators, and it has been derived to understand phenomena of collective synchronization in chemical and biological systems. Roughly speaking, the main features of this model are the following. Uncoupled oscillators run independently at their natural frequencies, and when the coupling is sufficiently weak, they still run incoherently.  At a critical value of the coupling strength, the system presents a phase transition to synchrony:  the oscillators spontaneously exhibit a collective behavior, that is
partial synchronization, the incoherent state loses stability and coherent dynamics emerge. Full synchronization occurs as the interaction strength goes to infinity. We refer to the review paper \cite{acerbon} for a detailed description of the model and for several related results. 
 
  
As the number of oscillators goes to infinity, the Mean Field approach comes into play. Recently, Carmona, Cormier and Soner \cite{carmonasoner} proposed a Mean Field Game (MFG) version of the classical Kuramoto model. The synergy between the Kuramoto and MFG formalisms has been already explored in \cite{carmona}, where a jet-lag recovery model was considered, and in \cite{yin}, where bifurcation arguments have been used to analyze incoherence and coordination in some large population game Kuramoto models. 
The main difference between classical and MFG Kuramoto models is the following; in the former, oscillators are treated as a particle system that evolve according to predetermined rules. In the latter, particles are rational agents, who are allowed to ``choose'' their evolution to minimize a (predetermined) cost depending on the evolution of the other oscillators. Equilibria are then considered in the Nash sense. 

A main question in these models is to understand the emergence of syncronization, and study its possible long time stability. 
In classical Kuramoto (Mean Field) models the evolution is naturally forward in time, and the long time analysis is quite well understood, se for instance \cite{Carrillo, Morales}. On the other hand, a main difficulty of the MFG setting is its forward-backward nature, because evolution runs forward while optimization runs backward by Bellman dynamic programming principle. 

In \cite{carmonasoner} the existence of a phase transition is observed, as in the classical Kuramoto model: first, for large interaction parameter, there are non-uniform stationary solutions, that become fully syncronized as the interaction parameter goes to $+\infty$. Moreover, the authors show that below a certain critical parameter, agents desynchronize: their distribution converges, in long time, to the uniform measure, at least for initial data in a neighborhood of the uniform distribution. In other words, the incoherent state (uniform distribution) is locally stable in long time when the interaction parameter is \textit{small} (or the discount factor is large). However, several interesting questions remain open. Is it possible to characterize all stationary equilibria? Is it possible to say something on their long time stability/instability?


In this paper we provide  some partial answers to the previous questions. In particular, we are able to describe stationary equilibria and study their local stability properties when the interaction parameter is \textit{large}.

\medskip

Let us now introduce the MFG version of the Kuramoto model,  in the periodic state space $\T=\R/2\pi\Z$, that will be identified with $(-\pi, \pi]$. The phase of a generic oscillator evolves according to \[dX_t= \alpha_t dt+\sqrt{2}dW_t,\]  
where $W_t$ is a standard Brownian motion. We first discuss the \textit{ergodic} framework (that is, when the discount parameter $\beta$ in \cite{carmonasoner} vanishes). In such case, the control $\alpha_t$ is chosen to minimize the long run cost
\begin{equation}\label{lran}\lim_{T\to+\infty}\frac{1}{T}\mathbb{E}\int_0^T \left[\frac{|\alpha_t|^2}{2}+2\kappa \int_{-\pi}^\pi \sin^2\left(\frac{X_t-y}{2}\right) dm(y)\right] dt,
\end{equation} where $m$ is the invariant measure of all the oscillators, that is, the observed stationary distribution of the environment, while $$\kappa > 0$$ is the interaction parameter. In an equilibrium regime, the law $\mathcal{L}(X_t^\alpha)$ of the generic oscillator, driven by the optimal control, converges as $t\to +\infty$ to the distribution $m$. Observe that 
\[ 2 \sin^2\left(\frac{x-y}{2}\right)=1-\cos(x-y)=1-\cos x\cos y-\sin x\sin y\] therefore the long run cost can be rewritten  as:
\begin{equation}\label{lran2}\lim_{T\to+\infty}\frac{1}{T}\mathbb{E}\int_0^T \left[\frac{|\alpha_t|^2}{2}-\kappa\cos X_t \int_{-\pi}^\pi \cos(y)dm(y)-\kappa\sin X_t \int_{-\pi}^\pi \sin(y)dm(y) \right] dt\, \, (+\kappa).
\end{equation}

Note that equilibria are translation invariant, that is: if $m$ is an equilibrium, then $m(\cdot-z)$ is also an equilibrium for every $z \in \R$. This is expected, because no syncronization to any ``special'' phase is enforced.

  
 Using the analytic (PDE) approach in MFG \cite{LL061, LL062, LL07}, the equilibrium regime in \eqref{lran2}  is encoded by $2\pi$-periodic solutions $(u,  \tilde \lambda, m)$ of the ergodic MFG
     \begin{equation}\label{mfgnonsim}
 \begin{cases} -u''+\frac{1}{2}|u'|^2+\tilde  \lambda= -\kappa\cos x  \int_{-\pi}^\pi \cos(y)m(y)dy -\kappa\sin x  \int_{-\pi}^\pi \sin(y)m(y)dy \quad \text{in $\T$}   \\ 
m(x) = \frac{e^{-u(x)}}{\int_{-\pi}^\pi e^{-u(y)}dy} \\
\text{$u,m$ are $2\pi$-periodic, $u(0) = 0$.}
\end{cases} 
 \end{equation} 
The density $m$ of the population of oscillators is a solution of the Fokker-Planck equation
\[
-m'' - (u'm)' = 0, \qquad \int_{-\pi}^\pi m = 1.
\]
It is well known, and easy to check, that the unique solution to the previous equation is explicitly given by the formula in \eqref{mfgnonsim}.
  
 We introduce the notion of incoherent and self-organizing solutions to \eqref{mfgnonsim}. 
\begin{definition}[Incoherent and self-organizing ergodic solutions]\upshape
The triple $(0,0, \frac{1}{2\pi})$ where $u\equiv 0\in \R$ is constant, $m\equiv \frac{1}{2\pi}$ is the uniform probability density on the torus is called the  \textit{incoherent} solution of the Kuramoto MFG \eqref{mfgnonsim}. 

A solution $(u,\lambda, m)$ to \eqref{mfgnonsim} is \textit{self-organizing} if it is not equal to the  incoherent  
solution. 
\end{definition} 
 
The first main result of the paper is the existence and \textit{uniqueness}, up to translation, of self-organizing solutions to \eqref{mfg}, for sufficiently large values of the interaction parameter $\kappa$. 
 
  \begin{theorem}\label{ex0}  There exists $\kappa_0>4$ such that, for all $\kappa\geq\kappa_0$, self organizing solutions to the Kuramoto system \eqref{mfgnonsim} are unique, up to translation in the $x$-variable.   \end{theorem}
 
Our uniqueness result answers positively, for $\kappa$ sufficiently large, to a conjecture proposed in \cite[Remark 7.4]{ carmonasoner}. We derive here some fine properties of a real valued function whose fixed points are connected with solutions of \eqref{mfgnonsim}. Note that such function is believed to be convex; here, we are able to obtain properties of its first derivative that are strong enough to classify all of its fixed points.

\bigskip

The second part of the paper is devoted to the study of the \textit{local} dynamical stability of self-organizing solutions to the Kuramoto MFG. 
Let us first briefly recall some known facts on the long time behavior of MFG. The classical Lasry-Lions monotone case is pretty well understood, see for instance \cite{Cannarsa, CPmaster, GMSR, PorMinMax} (and references therein), namely solutions enjoy an exponential turnpike property, that is: any solution $(u,m)=(u^T,m^T)$ of the finite $T$ horizon problem is exponentially close to the unique stationary state $(\bar u, \bar m)$ in the following sense:
\[
\|m(t) - \bar m\|_{L^2} + \|u(t) - \bar u\|_{L^2} \lesssim e^{-\omega t} + e^{-\omega (T-t)} \qquad \forall t\in[0,T].
\]
Such property is \textit{global}, namely it holds for solutions satisfying arbitrary initial-final condition. The core principle behind this kind of long-time stability is, in our viewpoint, the following. If the coupling is monotone and the Hamiltonian is uniformly convex, one can show that the quantity
\[
\Phi(t) = \|m(t) - \bar m\|^2_{L^2} + \|u(t) - \bar u\|^2_{L^2}
\]
satisfies the following inequality
\begin{equation}\label{Phi}
\int_{t_1}^{t_2}\Phi(t) dt \lesssim \Phi(t_1) + \Phi(t_2)
\end{equation}
for every $t_1 \le t_2 \in [0,T]$.  From this, it is possible to deduce the exponential decay: we detail the argument in the Appendix, Lemma \ref{expdecaylemma}. The inequality \eqref{Phi} is a straightforward consequence of the standard duality identity between \textit{state} $m$ and \textit{co-state} $u$, plus an application of the Poincar\'e inequality. It has then been noted in \cite{CiPo} that \eqref{Phi} is available also if the coupling is \textit{mildly} nonmonotone, at least in problems with nondegenerate diffusion. Indeed, the stabilization properties of these diffusions can be quantified in order to compensate a mild nonmonotone coupling. This observation will be important also in this work, as we will see below.

Still, long time stability in MFG is mainly understood whenever global uniqueness of dynamic and stationary equilibria holds, and for problems which are set on bounded domains. We are aware of a few exceptions only: \cite{BK} obtains stability for some deterministic problems with particular structure, and \cite{Gueant} studies the local stability for a special nonmonotone problem, for which a linear stability analysis can be carried out explicitly. In fact, a stable long time behavior in MFG that have no monotone structure is in general not expected \cite{CardaMaso, CCBrake, CJDE, CNur, Maso}.

Thus, the study of local stability of stationary solutions in MFG is   widely open when uniqueness fails. Local stability is actually what one would like to understand in a Kuramoto MFG, where, as we prove  in the first part of this paper, there is a continuum of different stationary solutions.

To simplify the stability analysis, we will restrict to \textit{even} solutions, and $\kappa$ large enough. Within this framework, we have shown that there exist \textit{two} stationary solutions only: the incoherent one and a self-organizing one (satisfying $\int \cos \bar m > 0$). The dynamic, finite-horizon version of the Kuramoto MFG in the time-space cylinder $(0,T)\times (-\pi, \pi)$ is: 
\begin{equation}\label{mfgev1}  \begin{cases} -u_t-u_{xx}+\frac{1}{2}|u_x|^2=  -\kappa\cos x  \int_{-\pi}^\pi \cos(y)m(t,y)dy
\\ m_t-m_{xx}-(mu_x)_x=0\\ m_x(t, \pi )=m_x(t, -\pi )=0\qquad u_x(t, \pi )=u_x(t, -\pi)=0\\ 
\text{$u(\cdot,t),m(\cdot,t)$ are even, $\int_{-\pi}^{\pi} m(x,t)dx=1, m(\cdot,t) \ge 0$ for all $t$}\end{cases}   \end{equation}

Since there can be several dynamic equilibria $(u,m)$ for any fixed initial-final conditions $m(0), u(T)$, our goal is to show the following local stability property / local (exponential) turnpike of the self-organizing solution $(\bar u, \bar m)$: there exists a neighborhood 
$\cU$ of $( \bar m, \bar u)$,  such that, for  any dynamic solution $(u,m)$ to \eqref{mfgev1} remaining in $\cU$ for all times, that is $(m(t),u(t)) \in \cU$ for all $t \in [0,T]$, it is true that   \begin{equation}\label{phi2}
\Phi(t) \lesssim e^{-\omega t} + e^{-\omega (T-t)} \qquad \forall t\in[0,T],
\end{equation}
for some constants $\omega$ independent of $T$. Here $\Phi(t)$ should be a positive function which quantifies the distance between $m(t)$ and $\bar m$ (and also between $u(t)$ and $\bar u$).

We are able here to identify a suitable $\Phi$; our second main result reads then informally as follows:
\begin{theorem*} Let $\kappa$ be large enough so that $(\bar u, \bar m)$ is the unique even self-organizing solution (Theorem \ref{ex0}). Assume that $(m,u)$ is a solution to \eqref{mfgev1} such that 
\[
m(x,t) \le C \bar m(x) \qquad \text{for all $x,t$.}
\]
Then, \eqref{phi2} holds with $\Phi(t) = \|m(t) - \bar m\|_{L^2(\bar m^{-1})}$.
\end{theorem*}
The precise results is stated in Theorem \ref{conv1}. Note that the constants involved in the estimate depend on $C, \kappa$, and $u_x, m|_{t = 0,T}$, but not on $T$. The result is obtained starting from the crucial observation that the rescaled variables $w(t,x)= u(t\kappa^{-\frac{1}{2}}, x\kappa^{-\frac14})$, $\mu(t,x)= \kappa^{-\frac14} m(t\kappa^{-\frac{1}{2}}, x\kappa^{-\frac14})$ solve a MFG
system where the coupling (formally) vanishes as $\kappa \to +\infty$, see \eqref{mfgresc}. Since the coupling is mild in this new scale, one is tempted to argue as in \cite{CiPo}, but soon realizes that a main difficulty is that the problem is set on a domain that becomes the real line in the limit $\kappa \to \infty$. We have then to implement weighted Poincar\'e inequalities, and stability of the Fokker-Planck equation in wighted $L^2$ spaces (see, for instance, \cite{BGG} and references therein). Though we address a specific problem, we believe that the functional setting used here can be useful to study the long time behavior in more general MFG which are set on unbounded domains (like the whole Euclidean space), for which, even in the Lasry-Lions monotone case, there are very few available results: we are only aware of \cite{Arapostathis}.

Note that an estimate like \eqref{phi2} just guarantees that trajectories that remain close to the equilibrium actually converge to it very quickly. Given an initial-final condition, the existence of these trajectories is not addressed here. Nevertheless, the kind of estimates that we obtain can be used to set up a topologic fixed point argument, which in fact yields existence, at least for boundary data that are close to the equilibrium, as in \cite{CiPo}.
Finally, the question of the long time local stability remains open at this stage for dynamic equilibria which are \textit{not even}. In this case, it is not clear whether or not they stabilize in long time to a stationary self-organizing one, and if so, which one of the infinitely many is selected. We believe that to tackle this issue one should have a look at \textit{orbital stability}, a key stability concept in Hamiltonian systems such as the Schr\"odinger equation. We will pursue this investigation in a future work.

%

\bigskip

The paper is organized as follows. In Section \ref{sezstat} we provide existence and uniqueness of symmetric self-organizing solutions to the problem for $\kappa$ larger than a threshold $\kappa_0>4$.  Section \ref{sezdinamica} contains the proof of the local stability of the self-organizing solutions. We finally collect in the appendix some useful estimates, and the proof of the Poincar\'e weighted inequality.  

 \subsection*{Acknowledgements}  The authors are members of GNAMPA-INdAM.  They were partially supported by the King Abdullah
University of Science and Technology (KAUST) project CRG2021-4674 ``Mean-Field Games: models, theory, and computational aspects''.
 
 \section{Ergodic self-organizing equilibria} \label{sezstat} 

In this section we prove Theorem \ref{ex0}. Most of the efforts will be devoted to classify \textit{even} solutions, that is, to show the following result.

\begin{theorem} \label{ex2}  There exists $\kappa_0>4$ such that for all $\kappa\geq\kappa_0$
 the Kuramoto MFG \eqref{mfgnonsim} admits, besides the incoherent solution, a unique \underline{even} self organizing solution $(u,\lambda, m)$ with $\int_{-\pi}^\pi m(x)\cos xdx>0$ and a unique \underline{even} self organizing solution $(u,\lambda, m)$ with $\int_{-\pi}^\pi m(x)\cos xdx<0$.   \end{theorem}
 
 Indeed, let $(u,\lambda, m)$ be any solution to \eqref{mfgnonsim}. Up to translation, we can always assume that
 \begin{enumerate}
 \item $\int_{-\pi}^\pi  m(y)\sin y \ dy=0$,
 \item $\int_{-\pi}^\pi m(y)\cos y \ dy \ge 0$,
 \item $u(x)=u(-x)$, $m(x) = m(-x)$ for all $x$.
 \end{enumerate}
To check (i), consider $\hat u(x) = u(x + z)-u(z), \hat m(x) = m(x + z)$, which still solves
  \[
 \begin{cases} -\hat u''+\frac{1}{2}|\hat u '|^2+\tilde  \lambda= -\kappa\cos x  \int_{-\pi}^\pi \cos(y)\hat m(y)dy -\kappa\sin x  \int_{-\pi}^\pi \sin(y)\hat m(y)dy \\ \qquad\qquad\qquad\qquad= 2\kappa \int_{-\pi}^\pi \sin^2\left(\frac{x-y}{2}\right) d\hat m(y) - \kappa  \\ 
\hat m(x) = \frac{e^{-\hat u(x)}}{\int_{-\pi}^\pi e^{-\hat u(y)}dy} \\
\text{$\hat u, \hat m$ are $2\pi$-periodic, $\hat u(0) = 0$.}
\end{cases}
\]
In addition,
\[
\int_{-\pi}^\pi \hat m(y)\sin y \ dy= \cos z \int_{-\pi}^\pi m(y)\sin y \ dy - \sin z \int_{-\pi}^\pi m(y)\cos y \ dy = 0
\]
for a suitable choice of $z$.

Regarding (ii),  if $\int_{-\pi}^\pi m(y)\cos y \ dy \le 0$ one can proceed as before by considering $\hat u(x) = u(x + \pi)-u(\pi), \hat m(x) = m(x + \pi)$, which solves the same problem and satisfies also
\[
\int_{-\pi}^\pi \hat m(y)\cos y \ dy= - \int_{-\pi}^\pi m(y)\cos y \ dy \ge 0.
\]

Finally, if $(i)$ holds, then $(iii)$ holds as well, that is, $u$ and $m$ are even. Indeed, $u$ turns out to be a $2\pi$-periodic solution of the ergodic HJ equation
\[ -u''+\frac{1}{2}|u'|^2+\tilde  \lambda= -\kappa\cos x  \int_{-\pi}^\pi \cos(y)m(y)dy \qquad \text{and $u(0)=0$}. \]
Periodic solutions of the previous equation are known to be unique, namely the couple $(\tilde  \lambda, u)$ is unique. Since $u(-x)$ also solves the previous problem, we get that $u(x) = u(-x)$. Hence $u$ is even, and $m$ needs to be even as well.
 
 Therefore, any solution to \eqref{mfgnonsim} satisfies, \textit{up to translation}, the following problem:
 \begin{equation}\label{mfg}
 \begin{cases} -u''+\frac{1}{2}|u'|^2+\tilde  \lambda= -\kappa\cos x  \int_{-\pi}^\pi \cos(y)m(y)dy   \\ 
m(x)=\frac{e^{-u(x)}}{\int_{-\pi}^\pi e^{-u(y)}dy}, \\
\text{$\hat u, \hat m$ are even, $2\pi$-periodic, $\hat u(0) = 0$, $\int_{-\pi}^\pi m(y)\cos y \ dy \ge 0$}.
\end{cases} 
 \end{equation}
 Theorem \ref{ex2} states that self-organizing solutions to the previous problem are unique, hence  Theorem \ref{ex0} follows as a straightforward consequence.
 
\bigskip

We now proceed with the proof of Theorem \ref{ex2}. First of all we slightly rewrite the system \eqref{mfg} in an equivalent way. 
Define \[V(x):=1-\cos x\] and $\lambda:=\tilde \lambda+\kappa \int_{-\pi}^\pi m(y)\cos ydy$. Then \eqref{mfg} becomes 
  \begin{equation}\label{mfgnew}
 \begin{cases} -u''+\frac{1}{2}|u'|^2+  \lambda=  \kappa V(x)\left[1-\int_{-\pi}^\pi V(y)m(y)dy\right]   \\ u'(\pm \pi)=0\\ 
m(x):=\frac{e^{-u(x)}}{\int_{-\pi}^\pi e^{-u(y)}dy}.
\end{cases} 
 \end{equation} 
Since $u,m$ are even, it will be indeed convenient below to work with Neumann boundary conditions at the boundary of the set $(-\pi, \pi)$.
We say that $(u, \lambda, m)$ is a solution to \eqref{mfgnew} if $\lambda\in \R$, $(u,m)$ are smooth and
solves in the classical sense the first equation in \eqref{mfgnew} (classical solution are in fact $C^\infty$). 
Note that
\begin{equation}\label{Vbound}
\frac{x^2}6 \le V(x) \le \frac{x^2}2 \qquad \text{on $[-\pi, \pi]$}.
\end{equation}

   \begin{remark}[The rescaled system]\label{rescaled} \upshape Several arguments below exploit a blow-up of \eqref{mfgnew}. Let $w(x)=u(x \kappa^{-\frac{1}{4}})$ and $\mu(x)=\kappa^{-\frac{1}{4}}m(x \kappa^{-\frac{1}{4}})$. The rescaled problem then reads:
\begin{equation}\label{hjbresc}  -w''+\frac{1}{2}|w'|^2+  \kappa^{-\frac{1}{2}}\lambda=  V_\kappa(x) \left[1-\kappa^{-\frac{1}{2}}\int_{-\pi\kappa^{\frac14}}^{\pi\kappa^{\frac14}} V_\kappa(y)\mu(y)dy\right]\qquad\text{with }
\mu(x)=\frac{e^{-w(x)}}{\int_{-\pi\kappa^{\frac14}}^{\pi\kappa^{\frac14}} e^{-w(y)}dy}
\end{equation} where $V_\kappa(x)=\kappa^{\frac{1}{2}}V(x \kappa^{-\frac{1}{4}})$. 
Note that by \eqref{Vbound}, $V_\kappa$ satisfies
\begin{equation}\label{Vkbound}
\frac{x^2}6 \le V_k (x) \le \frac{x^2}2 \qquad \text{on $(-\pi \kappa^{\frac{1}{4}}, \pi \kappa^{\frac{1}{4}})$}.
\end{equation}
The main advantage of this blow-up is that it ``weakens'' the coupling between $u$ and $m$, since $\kappa$ in front of $\int m \cos$ becomes $\kappa^{-\frac{1}{2}}$.
\end{remark}

In order to prove uniqueness of solutions to \eqref{mfgnew} (and in fact also existence), we note that such solutions correspond to fixed points of a function of a real variable. Fix $a\in \R$, and consider the solution $(u_a, \lambda_a)$  to the ergodic Hamilton-Jacobi equation with Neumann boundary conditions
\begin{equation}\label{hjb}  -u_a''+\frac{1}{2}|u_a'|^2+  \lambda_a=  \kappa V(x)\left[1-a\right], \qquad u_a(0)=0
\end{equation}  where the parameter $\kappa$ is fixed, 
and define $F_\kappa:\R \to \R$ as 
\begin{equation}\label{function}F_\kappa(a):= \int_{-\pi}^\pi V(y)m_a(y)dy  \qquad \text{ where }m_a(y):=\frac{e^{-u_a(y)}}{ \int_{-\pi}^\pi  e^{-u_a(y)}dy}. \end{equation}
It is well known, see \cite{cdm, c,ll} that   for every $\kappa, a\in \R$ there exists a unique  $\lambda_a\in \R$  and a smooth (even) function $u_a$ which solves in the classical sense \eqref{hjb}. 

Note that
\[
F_k(a) = 1 - \int_{-\pi}^\pi m_a(y)\cos y \ dy \in [0,2],
\]
hence $F_k : [0,2] \to [0,2]$. The stability property of the ergodic problem with respect to variations of the parameters is a well known result, see e.g. \cite[Proposition 3]{cdm}, hence $F_k$ is also continuous. Clearly, there is a one-to-one correspondence between fixed points $a = F_k(a)$ and solutions to \eqref{mfgnew}. Note that since we consider solutions $\int_{-\pi}^\pi m_a(y)\cos y \ dy \ge 0$ we should restrict to
\[
a = 1 - \int_{-\pi}^\pi m_a(y)\cos y \ dy \in [0,1].
\]
We shall see below that actually $F_k : [0,1] \to [0,1]$. In particular $a=1$ is a fixed point  of $F_\kappa$ for every $\kappa$, and $(u_1, \lambda_1, m_1)$ coincides with the incoherent solution to the Kuramoto MFG.

We now derive a crucial representation formula for $F'_\kappa$.

%
%
%

\begin{proposition}[Variations with respect to $a$]\label{propositionderivative}
 
Let $a\in [0,2]$ and $(u_a, \lambda_a)$ the solution to \eqref{hjb} with parameter $\kappa$. 

Then there exists a unique $\lambda_a' \in \R$ and a unique smooth $v_a$ which solves in classical sense the equation 
\begin{equation}\label{hjb1}  -v_a''+v'_au'_a+  \lambda'_a=  -\kappa V(x), \quad v_a(0)=0  \end{equation}
with Neumann boundary conditions. Moreover 
\begin{equation}\label{derivata1}\lambda'_a=\lim_{h\to 0}\frac{\lambda_{a+h}-\lambda_a}{h} =-\kappa F_\kappa(a)= - \kappa \int_{-\pi}^\pi  V(y)m_a(y)dy\end{equation}
and
\begin{eqnarray}  F_\kappa'(a)= \lim_{h\to 0}\frac{F_\kappa(a+h)-F_\kappa(a)}{h}&=& -\frac{1}{\kappa} \lim_{h\to 0}\frac{\lambda'_{a+h}-\lambda'_a}{h}\nonumber   \\ \label{derivata2}   & =&   \frac{1}{\kappa} \int_{-\pi}^\pi (v_a'(y))^2m_a(y)dy
\\ \label{derivata22} & =&  - \int_{-\pi}^\pi \left( \frac{\lambda_a' }{\kappa}+V(y)\right) v_am_a dy.\end{eqnarray}
 \end{proposition}
 
In particular $F_\kappa$ is a nondecreasing function, and since $F_k(1)=1$, then $F_\kappa : [0,1] \to [0,1]$. Observe also that for $a=1$, since $u_1=0, \lambda_1=0$, $m_1=\frac{1}{2\pi}$, we get $$\lambda'_1 = -\kappa, \qquad v_1(x)=\kappa\cos x, \qquad F_\kappa'(1)=\frac\kappa 2.$$

 \begin{proof}  To obtain \eqref{hjb1} and \eqref{derivata1} we follow the same arguments as in \cite[Section 5.1]{cdm}. First of all, the existence of a unique couple $(\lambda', v)$ solving \eqref{hjb1} is standard.
 
  Fix $h$ small and consider $(u_{a+h}, \lambda_{a+h})$ the solution to \eqref{hjb} with parameter $\kappa$. Define $v_h=\frac{u_{a+h}-u_a}{h}$ and $ \lambda_h=\frac{\lambda_{a+h}-\lambda_a}{h}$.  Then $v_h, \lambda_h$ solves
\[ -v''_h+\frac{1}{2}(u_a'+u_{a+h}')v'_h +\lambda_h=-\kappa V(x). \] If we multiply by $m_a$ this equation and integrate, recalling that $-m_a''-(m_au_a')'=0$, we get 
\[\lambda_h=-\kappa \int_{-\pi}^\pi V(y)m_a(y)-\frac{1}{2} \int_{-\pi}^\pi h (v_h')^2 m_ady.\]
Moreover $m_h=\frac{m_{a+h}-m_a}{h}$ solves 
\[-m_h''-(u_{a+h}' m_h+v_h'm_a)'=0\] with $\int_{-\pi}^\pi m_h=0$. If we multiply by $m_h$ the equation for $v_h$ and we subtract to it the equation for $m_h$ multiplied by $v_h$,   and we  integrate over $(-\pi, \pi)$ we conclude that 
\[\frac{1}{2} \int_{-\pi}^\pi h (v_h')^2(m_a+m_{a+h})dy=\kappa  \int_{-\pi}^\pi  V(y)(m_{a+h}-m_a)dy.\]
So, since $m_{a+h}-m_a\to 0$ uniformly as $h\to 0$,  we conclude that $\lim_{h\to 0}\lambda_h= -\kappa \int_{-\pi}^\pi V(y)m_a(y)=-\kappa F_\kappa(a)$. Moreover since $\frac{1}{2}(u_a'+u_{a+h}')\to u_a'$ as $h\to 0$ uniformly, 
by stability of viscosity solutions and uniqueness of the ergodic constant we conclude that  $\lambda'=\lim_{h\to 0} \lambda_h$ and that $\lim_h v_h=v_a$ uniformly in $C^{1,\gamma}$. 

To obtain \eqref{derivata2} and \eqref{derivata22} we now define  $z_h=\frac{v_{a+h}-v_a}{h}$ and $ \lambda'_h=\frac{\lambda'_{a+h}-\lambda'_a}{h}$.  Then $z_h, \lambda'_h$ solves
\[ -z''_h+ u_{a+h}'z_h'+\frac{u_{a+h}'-u_{a}'}{h} v_a' + \lambda'_h=0. \]  
We multiply this equation by $m_{a+h}$ and integrate in $[-\pi,\pi]$ recalling that $-m_{a+h}''-(m_{a+h}u_{a+h}')'=0$ and we obtain:
\[- \lambda'_h=\int_{-\pi}^\pi  \frac{u_{a+h}'-u_{a}'}{h} v_a' m_{a+h}  dy=\int_{-\pi}^\pi  v_h' v_a' m_{a+h}   dy. \] As $h\to 0$, we get that $v_h\to v_a$ in $C^{1,\gamma}$ and $m_{a+h}\to m_a$ uniformly, so, passing to the limit we obtain
\begin{equation}\label{lambdasec}
\big(-\kappa F_\kappa(a)\big)'=\lim_h\lambda_h'= - \int_{-\pi}^\pi (v_a'(y))^2m_a(y)dy.
 \end{equation}

Finally, we multiply \eqref{hjb1} by $m_a v_a$ and subtract the equation $-m_a''-(m_au_a')'=0$ multiplied by $\frac{1}{2}v_a^2$ and integrate:
\begin{eqnarray*}
&&\int_{-\pi}^\pi  \left[-v_a''+v'_au'_a+  \lambda'_a+\kappa V(x) \right] v_a m_a-\left[-m_a''-(m_au_a')'\right]\frac{v_a^2}{2}dy\\
&=& \int_{-\pi}^\pi \left[-\frac{1}{2}(v_a^2)''m_a+(v_a')^2m_a+ v'_au'_a v_a m_a+\left[\lambda'_a+\kappa V(x) \right] v_a m_a+\frac{1}{2} (v_a^2)''m_a-m_au_a'v_av_a'\right] dy\\
&=& \int_{-\pi}^\pi \left[ (v_a')^2m_a+\left[\lambda'_a+\kappa V(x) \right] v_a m_a\right] dy=0.
\end{eqnarray*}
This, together with \eqref{lambdasec}, implies \eqref{derivata22}. 
 \end{proof} 
%

\begin{remark}[Symmetry of $F_\kappa$]\upshape\label{amaggiore1} Let us observe that if $u_a, \lambda_a$ is the solution to \eqref{hjb} associated to $a\in [0,1]$ (where we fix $u_a(0)=0$), then
$\tilde u_a(x):=u_a(x+\pi)-u_a(\pi), \tilde \lambda_a:=\lambda_a-2\kappa(1-a)$  is the  solution to \eqref{hjb} associated to $\tilde a=2-a$. 
This implies that for all $a\in [0,2]$, 
\[u_{2-a}(x):=u_a(x+\pi)-u_a(\pi)\qquad \lambda_{2-a}=\lambda_a-2\kappa(1-a)\qquad F_\kappa(2-a)=2-F_\kappa(a).\]
\end{remark}

\begin{figure}
\centering
\includegraphics[width=4cm]{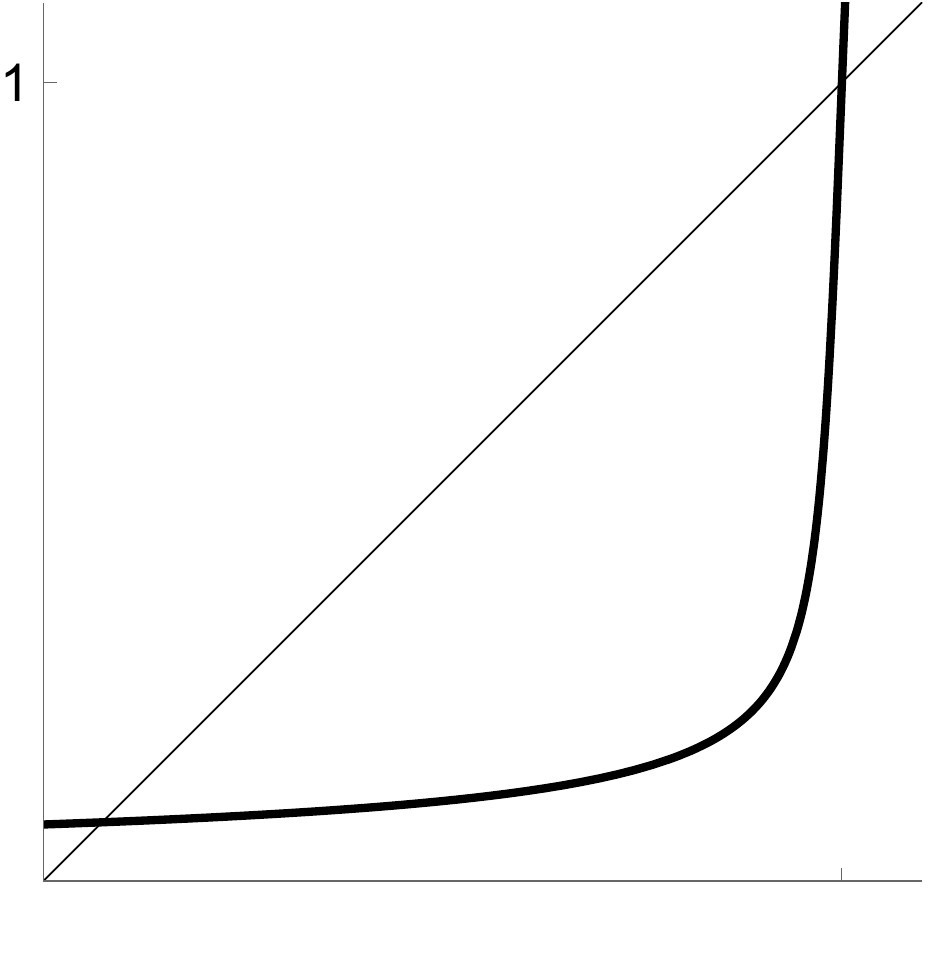}
\caption{\footnotesize Plot of the function $a \mapsto F_k(a)$, for large $\kappa$.}\label{Fk}
\end{figure}

 We now are going to show that if $\kappa$ is sufficiently large then there exists a fixed point $\bar a\in (0,1)$ of the function $F_\kappa$, and that this fixed point is unique. This would also imply, by the previous Remark \ref{amaggiore1}, that for  $\kappa$ is sufficiently large, $2-\bar a$ is the unique fixed point of the function $F_\kappa$ in $(1,2)$. Our strategy is to derive the following properties of $F_\kappa$:
 \begin{itemize}
 \item Theorem \ref{ex1}: for any $\delta > 0$, $F'_\kappa = O(\kappa^{-1/2})$  on $[0,1-\delta]$ for large enough $\kappa$.
 \item Proposition \ref{propositiontau}: there exists $\tau_0 \in (0,1)$ such that $F_\kappa'\left(1-\frac{\tau}{\kappa}\right)\geq \frac{\kappa}{4}$ for any $\tau \in [\tau_0, 1]$ and $\kappa$ large enough.
 \end{itemize}
 The two points above show that there can be just one fixed point in an interval $[0, 1-\delta]$ and close to $a=1$. In other words, for large $\kappa$, $F_\kappa$ is almost flat (and close to zero) for all $a$, and as soon as $a$ gets close to one, $F_\kappa$ abruptly reaches the fixed point $a=1$, see Figure \ref{Fk}. Combining this information with the fact that $F_\kappa$ is monotone gives the result.

Below, $a\in [0,1]$, and $(u_a, \lambda_a)$ the solution to \eqref{hjb} with $u_a(0)=0$ and with interaction parameter $\kappa \ge 1$. Most importantly, positive constants in the statements will be independent of $\kappa$.

\begin{proposition}\label{propositionstimal} There exists $\ell > 0$ such that
\begin{equation}\label{stimalambda} 0\leq \lambda_a  \leq \ell \kappa^{\frac{1}{2}}. \end{equation} 
\end{proposition} 

\begin{proof} Let $\varphi(x):=e^{-\frac{u_a(x)}{2}}$. Then, the couple $(\varphi, \lambda_a)$ solves
\[
\begin{cases}
-\varphi''(x) + \frac{\kappa  V(x) [1-a]}2 \varphi(x) = \lambda_a \varphi(x) & \text{on $(-\pi, \pi)$}, \\
\varphi'(\pi) = \varphi'(-\pi) = 0.
\end{cases}
\]
Since $\varphi > 0$ on $(-\pi, \pi)$, $\lambda_a$ is the first (nontrivial) eigenvalue of the Schr\"odinger operator $-\Delta + \frac{\kappa V(x) [1-a]}2$ on $(-\pi, \pi)$ (with Neumann boundary conditions), hence it has the  following well-known characterization
\[
\lambda_a = \inf_{\substack{\phi \in H^1(-\pi, \pi) \\ \int_{-\pi}^\pi \phi^2 = 1}} \int_{-\pi}^{\pi}  |\phi'|^2 + \frac{\kappa V(x) [1-a]}2 \phi^2 dx,
\]
which yields $\lambda_a \ge 0$ as a straightforward consequence. Pick any smooth nonnegative $\psi$ with compact support in $(0,1)$ and such that $\int_{-\pi}^\pi \psi^2 = 1$, and let
\[
\phi(x) = \kappa^{1/8} \psi\left(x \kappa^{1/4} \right).
\]
Clearly, $\phi$ has support in $(0,\kappa^{-1/4})$ and it satisfies $\int_{-\pi}^\pi \phi^2 = 1$. Therefore,
\[
\lambda_a \le \int_{0}^{\kappa^{-1/4}} | \phi'|^2 + \frac{\kappa V(x) [1-a]}2 \phi^2 dx= \int_{0}^{1} \kappa^{1/2} |\psi'|^2 + \frac{\kappa V(x \kappa^{-1/4}) [1-a]}2 \psi^2dx.
\]
Using \eqref{Vbound} we obtain
\[
\lambda_a \le \kappa^{1/2} \int_{0}^{1}  | \psi'|^2 + \frac{x^2}4 \psi^2dx,
\]
that gives the conclusion.
\end{proof}

\begin{remark} [Mathieu functions] \normalfont As one can see from the previous proof, $\phi = e^{-u_a}$ is a so-called \textit{Mathieu function}, because it solves an equation of the form
\[
-\varphi'' + (b+q \cos x) \varphi = 0
\]
for some real $b,q$; $b=b(q)$ is the characteristic number, and it is strictly related to $\lambda_a$ in our formulation (while $q$ is proportional to $a$ and $\kappa$). This class of special functions has been extensively studied during the last century \cite{mathieu}. For instance, by known results one could infer very precise asymptotics of $\lambda_a$ as $\kappa \to \infty$. One could then prove that $F_\kappa$ has just two fixed points on $[0,1]$ by showing for instance that it is convex (which is reasonable if one looks at Figure \ref{Fk}). Since $F_\kappa'' = - \lambda_a'''/k$ by Proposition \ref{propositionderivative}, this amounts to establish the sign of the third derivative of $\lambda_a$ with respect to $a$. Unfortunately, we are not aware of any result on the behavior of derivatives of  the characteristic number $b$ as a function of $q$.

\end{remark}

\begin{proposition}\label{uabobelo} Let $\delta\in (0,1)$ and $a\in [0,1-\delta]$. 
Then, there exists $\tilde\kappa=\tilde\kappa(\delta)$ (with $\tilde \kappa(\delta)\to +\infty$ as $\delta\to 0$)  such that if $\kappa\geq \tilde\kappa(\delta)$ then 
\[
c_1 \kappa^{1/2} x^2 - c_3 \le u_a(x) \le  c_2 \kappa^{1/2} x^2 + c_3  \qquad \text{for $|x| \le \pi$}.
\]  for some $c_2, c_3>0$ and $c_1=c_1(\delta)>0$ with $c_1(\delta)\to 0$ as $\delta\to 0$. 
\end{proposition} 

\begin{proof}
We first need to obtain a control on $u_a$ close to $x=0$. We rescale the equation \eqref{hjb} as  in Remark \ref{rescaled}: we let  $w_a(x)=u_a(x \kappa^{-\frac{1}{4}})$ and we get 
\begin{equation}\label{hjbresc1}  -w_a''+\frac{1}{2}|w_a'|^2+  \kappa^{-\frac{1}{2}}\lambda_a=  \kappa^{\frac{1}{2}}V(x \kappa^{-\frac{1}{4}}) [1-a].
\end{equation} Observe that $ \kappa^{\frac{1}{2}}V(x \kappa^{-\frac{1}{4}})$ satisfies \eqref{Vkbound}, hence it is locally bounded with respect to $x$, uniformly in $\kappa$. 
By the local gradient estimates for $w_a$ (see e.g. \cite{c, ll}), we have that $|w_a'|\leq C_r$ on any interval $[-r, r] \subset (-\pi \kappa^{1/4}, \pi \kappa^{1/4})$, which implies, going back to $u_a$, that  $|u_a'|\leq C_r\kappa^{\frac{1}{4}}$ on $[-r \kappa^{-1/4}, r \kappa^{-1/4}]$. Hence,
\begin{equation}\label{uabr}
|u_a(x)| \le r C_r \qquad \text{on $[-r \kappa^{-1/4}, r \kappa^{-1/4}]$}.
\end{equation}

\smallskip

Let us now proceed with the bound for $u_a$ from above, by constructing a suitable supersolution of the HJ equation. Note that $u_a$ is even, hence it suffices to argue on $[0, \pi]$.
Let
\[
\bar u(x) := \frac{\bar c}2 \kappa^{1/2} x^2, 
\]
where $\bar c \ge 1$ will be chosen below (large). We have that, for $x \ge r \kappa^{-1/4}$, 
\begin{multline*}
-\bar u''+\frac{1}{2}|\bar u'|^2+  \lambda_a -  \kappa V(x)\left[1-a\right] \stackrel{\lambda_a \ge 0}{\ge} 
-\bar c \kappa^{1/2} + \frac{\bar c^2 \kappa}2 x^2 - \kappa V(x)\left[1-a\right]  \stackrel{\eqref{Vbound}}{\ge} \\
-\bar c \kappa^{1/2} + \frac{\kappa}2(\bar c^2 - 1) x^2 \ge
\kappa^{1/2} \big(-\bar c + (\bar c^2-1)r^2\big) \ge 0,
\end{multline*}
provided that $\bar c = \bar c(r)$ is chosen large enough. Note that $\bar u'(\pi) > 0$. By the Maximum Principle, the maximum of $u_a-\bar u$ on $r \kappa^{-1/4} \le x \le \pi$ is achieved at the boundary. Since it cannot be achieved at $x=\pi$ (that would contradict Hopf's Lemma), we get recalling  \eqref{uabr} that
\[
u_a(x)-\bar u(x) \le u_a(r \kappa^{-1/4}) -\bar u(r \kappa^{-1/4}) \le rC_r \qquad \text{in $r \kappa^{-1/4} \le x \le \pi$},
\]
which, using again \eqref{uabr}, yields
\[
u_a(x) \le  \frac{\bar c}2 \kappa^{1/2} x^2 + rC_r  \qquad \text{on $|x| \le \pi$}.
\]
Pick now any $r$ ($r=1$ for instance) to conclude the bound on $u_a$ from above.

\smallskip

To control $u_\alpha$ from below, we use the following subsolution on $x \ge r \kappa^{-1/4}$:
\[
\underline u(x) := \frac{c}2 \kappa^{1/2} (1-\cos x),
\]
for $c > 0$ small. Note that now $r$ will need to be chosen large enough so that $\ell - \frac{\delta r^2}{12} \le -1$. For $x \ge r \kappa^{-1/4}$, and $c$ small so that $c^2 \le \delta / 2$ and $c \le 1$,
\begin{multline*}
-\underline u''+\frac{1}{2}|\underline u'|^2+  \lambda_a -  \kappa V(x)\left[1-a\right] \stackrel{\eqref{stimalambda}}{\le} \\
c \kappa^{1/2}\cos x + \frac{c^2 k}2(1-\cos x)(1+\cos x) + \ell \kappa^{1/2} -  \kappa (1- \cos x)\left[1-a\right]  \le \\
c \kappa^{1/2} + \ell \kappa^{1/2} + \kappa (1-\cos x)\left(\frac{c^2}2(1+\cos x) - 1 + a\right) \le \\
c \kappa^{1/2} + \ell \kappa^{1/2} - \frac {\delta \kappa} 2  (1-\cos x) \le \kappa^{1/2} \left(c + \ell - \frac{\delta r^2}{12}\right) \le 0.
\end{multline*}
Arguing as before, we obtain
\[
u_a(x) \ge \underline u(x) \ge \frac{c \kappa^{1/2}}{12}  x^2- rC_r  \qquad \text{on $|x| \le \pi$}.
\]

  \end{proof} 
  
  \begin{corollary}\label{corollary2.9} Let $\delta\in (0,1)$. Fix $a\in [0,1-\delta]$. Then, for $\kappa\geq \tilde\kappa(\delta)$ (where $\tilde\kappa(\delta)$ is as in Proposition \ref{uabobelo}) there holds 
  \begin{equation}\label{mbound}
C^{-1}\kappa^{1/4}e^{-c_2 \kappa^{1/2} x^2}\leq  m_a(x) = \frac{e^{-u_a(x)}}{ \int_{-\pi}^\pi  e^{-u_a(x)}dx} \le C \kappa^{1/4 }e^{-c_1 \kappa^{1/2} x^2} \qquad \text{for all $x \in [-\pi, \pi]$},
  \end{equation} for some $C>0$ and for $c_1,c_2 $ as in Proposition \ref{uabobelo}. 
Moreover 
  \begin{equation}\label{Fbound}
    F_\kappa(a) = \int_{-\pi}^\pi V(x)m_a(x)dx \le \frac{C'}{\kappa^{1/2}}
  \end{equation}
  for some $C' >0$. Here, $C,C'$ depend on $c_1, c_2, c_3$.
  
  \end{corollary}
  
  \begin{proof} For the first assertion, it is sufficient to use Proposition \ref{uabobelo}:
  \[
  \frac{e^{-u_a(x)}}{ \int_{-\pi}^\pi  e^{-u_a(x)}dx} \le e^{2c_3}   \frac{e^{-c_1 \kappa^{1/2} x^2}}{\int_{-\pi}^\pi  e^{-c_2 \kappa^{1/2} x^2}dx} = e^{2c_3}   \kappa^{1/4 }\frac{e^{-c_1 \kappa^{1/2} x^2}}{\int_{-\pi \kappa^{1/4}}^{\pi \kappa^{1/4}}  e^{-c_2 y^2}dy}
  \]
    \[
  \frac{e^{-u_a(x)}}{ \int_{-\pi}^\pi  e^{-u_a(x)}dx} \ge e^{2c_3}   \frac{e^{-c_2 \kappa^{1/2} x^2}}{\int_{-\pi}^\pi  e^{-c_1 \kappa^{1/2} x^2}dx} = e^{2c_3}   \kappa^{1/4 }\frac{e^{-c_2\kappa^{1/2} x^2}}{\int_{-\pi \kappa^{1/4}}^{\pi \kappa^{1/4}}  e^{-c_1 y^2}dy}.
  \]
  To get the second one, since $V(x) \le x^2/2$,
  \[
  F_\kappa(a) \le C \int_{-\pi}^\pi x^2  e^{-c_1 \kappa^{1/2} x^2} \, \kappa^{1/4 } dx = \frac{C}{\kappa^{1/2}} \int_{-\pi \kappa^{1/4}}^{\pi \kappa^{1/4}} y^2  e^{-c_1 y^2} \, dy. 
  \]
  \end{proof}
  
  \begin{proposition}\label{vabobelo} Let $\delta\in (0,1)$. Fix $a\in [0,1-\delta]$, and consider   $(v_a, \lambda'_a)$ the solution to \eqref{hjb1} with $v_a(0)=0$ and with interaction parameter $\kappa$. Then,  there exists $\bar \kappa(\delta)\geq \tilde \kappa(\delta) $  (where $\tilde \kappa(\delta)$ is as in Proposition \ref{uabobelo})  such that for $\kappa\geq \bar\kappa(\delta)$   there holds 
  \begin{equation}\label{vabo}
|v_a(x)| \le  \bar c_1 \kappa^{1/2} x^2 + \bar c_2  \qquad \text{for $|x| \le \pi$}
\end{equation}
for some $\bar c_1, \bar c_2>0$. 
\end{proposition} 

\begin{proof} We start with some bounds on $v_a$ close to $x=0$. Since $v_a$ is a solution to \eqref{hjb1},  and $u_a$ is even, we get that also $v_a$ is even, and then $v_a'(0)=0$. By direct integration of \eqref{hjb1} we get  
\[
v'_a(x) = e^{u_a(x)} \int_0^x e^{-u_a(s)}\big(\lambda'_a + kV(s) \big) ds.
\]
Note that $|\lambda'_a| = \kappa F_\kappa (a) \le C\kappa^{1/2}$ by \eqref{Fbound}. Hence, by the control on $u_a$ obtained in Proposition \ref{uabobelo} we get, for any $r \ge 0$ and $|x| \le r \kappa^{-1/4}$,
\begin{multline*}
|v'_a(x)| \le e^{c_2\kappa^{1/2} x^2} e^{2c_3}\int_0^x \big(|\lambda'_a| + \kappa V(s) \big) ds \le e^{c_2 r^2} e^{2c_3} \int_0^x \left(C\kappa^{1/2} + \frac \kappa 2 s^2  \right)  ds \le\\  e^{c_2 r^2}  e^{2c_3} r \kappa^{-1/4} \left(C\kappa^{1/2} + \frac{\kappa^{1/2}} 2 r^2  \right) = C_r \kappa^{1/4},
\end{multline*}
which in turn yields
\begin{equation}\label{vabr}
|v_a(x)| \le r C_r \qquad \text{on $[-r \kappa^{-1/4}, r \kappa^{-1/4}]$}.
\end{equation}

Now we need to control from above and below $v_a$ in the annulus $r \kappa^{-1/4} \le |x| \le \pi$ by constructing suitable sub/supersolutions of \eqref{hjb1}. We first pick $r$ such that
\[
\frac{r^2 \delta}{12} \ge \ell, \qquad \frac{r^2}{6} \ge C+\ell,
\]
where $\ell, C$ appear in \eqref{stimalambda} and \eqref{Fbound} respectively. Let $\underline v(x) = -\eta u(x)$, where $\eta$ satisfies
\[
\eta \ge \frac2{\delta}.
\]
For $x \ge r \kappa^{-1/4}$ (note that $v_a$ is even, hence all the arguments below adapt to $x \le -r \kappa^{-1/4}$),
\begin{multline*}
-\underline v''+\underline v' u'_a+  \lambda'_a + \kappa V = \eta u_a''  - \eta |u_a'|^2 + \lambda'_a + \kappa V \stackrel{\lambda' \le 0}{\le }
\eta \lambda_a - \eta \kappa V[1-a]  + \kappa V \stackrel{\eqref{stimalambda}}{\le } \\
\eta \ell \kappa^{1/2} +\kappa V[1-\eta(1-a) ] \le \eta \ell \kappa^{1/2} - \kappa V \frac{\eta \delta}2 \le \eta \ell \kappa^{1/2} - \kappa x^2 \frac{\eta \delta}{12} \le
\eta \kappa^{1/2}\left(\ell - \frac{r^2 \delta}{12}\right) \le 0.
\end{multline*}
Hence $\underline v - v_a$ is a subsolution of \eqref{hjb1} on $x \ge r \kappa^{-1/4}$. By Hopf's Lemma (recall that $\underline v'(\pi) - v_a'(\pi) = 0$) and the maximum principle,
\[
\underline v(x) - v_a(x) \le \underline v(r \kappa^{-1/4}) - v_a(r \kappa^{-1/4}) \le rC_r \qquad \text{for $r \kappa^{-1/4} \le x \le \pi$},
\]
which implies
\[
v_a(x) \ge -\eta u(x) -rC_r \ge - \eta c_2 \kappa^{1/2} x^2 - \eta c_3-rC_r \qquad \text{for $|x| \le \pi$}
\]
in view of Proposition \eqref{uabobelo}.

To control $v_a$ from above, we argue similarly with $\bar v = \eta u$. Indeed, for $x \ge r \kappa^{-1/4}$,
\begin{multline*}
-\bar v''+\bar v' u'_a+  \lambda'_a + \kappa V = - \eta u_a''  + \eta |u_a'|^2 + \lambda'_a + \kappa V \stackrel{\eqref{Fbound}}{\ge }
- \eta \lambda_a + \eta \kappa V[1-a] - \eta  C \kappa^{1/2}  + \kappa V \ge  \\
 - \eta  (C+\ell) \kappa^{1/2}  + \kappa V \ge  - \eta  \kappa^{1/2} (C+\ell)  +  \kappa^{1/2} \frac{r^2}6
= \eta  \kappa^{1/2}\left(-C-\ell + \frac{r^2}6\right) \ge 0,
\end{multline*}
and we conclude as above that
\[
v_a(x) \le \eta u(x) + rC_r \le  \eta c_2 \kappa^{1/2} x^2 + \eta c_3-rC_r \qquad \text{for $|x| \le \pi$}.
\]

\end{proof}
Our first result is the existence and uniqueness of a fixed point of $F_\kappa$ on the set $[0,1-\delta$, for  $\kappa$ sufficiently large. 
  \begin{theorem}\label{ex1} 
  Let $\delta\in (0,1)$. Fix $a\in [0,1-\delta]$, and consider   the map $F_\kappa(a)$ defined in \eqref{function}, with interaction parameter $\kappa $. Then,  for $\kappa\geq \bar  \kappa(\delta)$ (where $\bar \kappa(\delta)$ is as in   Proposition \ref{vabobelo}) there holds  for some $C>0$ 
\[
0\leq F'_\kappa(a)\le C \kappa^{-\frac12}.\] 
In particular  there exists $\kappa_0(\delta)\geq \bar \kappa(\delta)$ (with $\kappa_0(\delta)\to +\infty$ when $\delta\to 0$), such that if $\kappa\geq \kappa_0(\delta)$, then 
$F_\kappa:[0,1-\delta ]\to [0,1-\delta]$ is  a contraction. Hence it admits a unique fixed point, which is associated to a self-organizing solution to the Kuramoto MFG \eqref{mfgnew}.
  \end{theorem}
  
 Note that, by  Remark \ref{amaggiore1}, $F_\kappa:[1+\delta,2 ]\to [1+\delta,2]$  is also a contraction.
  
\begin{proof}  

We recall  that $|\lambda'_a| = \kappa F_\kappa (a) \le C\kappa^{\frac12}$ by \eqref{Fbound} and $V(x) \le x^2 /2$. So
 \begin{eqnarray*}0\leq F'_\kappa(a) &\stackrel{\eqref{derivata22}}{=}& -\int_{-\pi}^\pi  \frac{\lambda_a' }{\kappa}v_am_a+ V(x) v_am_a dx  \\
&\leq & \int_{-\pi}^\pi  \left|\frac{\lambda_a' }{\kappa}+V(x)\right| |v_a| m_a  dx
\\ &\stackrel{\eqref{mbound}, \eqref{vabo}}{\le} &  \int_{-\pi}^\pi \left( \frac{C }{\kappa^{\frac{1}{2}}} +\frac{x^2}{2}\right)   (\bar c_1 \kappa^{\frac12} x^2 + \bar c_2)C' \kappa^{\frac14 }e^{-c_1 \kappa^{\frac12} x^2}   dx\\
&=&  \int_{-\pi\kappa^{\frac14}}^{\pi\kappa^{\frac14}}\left( \frac{C }{\kappa^{\frac{1}{2}}} +\frac{y^2}{2 \kappa^{\frac12}}\right)  ( \bar c_1 y^2 + \bar c_2)C' e^{-c_1 y^2}   dy \\
&\le&  \frac{C'' }{\kappa^{\frac{1}{2}}}   \int_{-\infty}^{\infty} (1+ y^2)( \bar c_1 y^2 + \bar c_2) e^{-c_1 y^2}   dy
\leq  \frac{C'''}{\kappa^{\frac12}}
\end{eqnarray*} 
for all $a \in [0,1-\delta]$. Now it is sufficient to choose $\kappa_1(\delta)\geq \bar \kappa(\delta)$ sufficiently large such that $C'''\kappa_1^{-1/2}(\delta)<1$. So, the map $F_\kappa$ is a contraction in $[0,1-\delta]$ and we conclude by Banach-Caccioppoli theorem. 
\end{proof} 

To conclude the proof of Theorem \ref{ex2}, we now show that the incoherent solution to the Kuramoto system is isolated as $\kappa$ is sufficiently large. 
\begin{proposition} \label{propositiontau} There exists $\tau_0\in (0,1)$ such that for all $\kappa>4$, there holds  for $\tau\in (0, \tau_0]$, 
\[ F_\kappa'\left(1-\frac{\tau}{\kappa}\right)\geq \frac{\kappa}{4}, \qquad \text{so that }\qquad F_\kappa\left(1-\frac{\tau}{\kappa}\right)\leq  1-\frac{\tau}{4}<1-\frac{\tau}{\kappa}.\]
Consequently, $F_\kappa$ in $\left[1-\frac{\tau_0}{\kappa},1+\frac{\tau_0}{\kappa}\right]$ admits a unique fixed point which is  $a=1$. 
\end{proposition} 
\begin{proof} 
We show that $F_\kappa$ in $\left[1-\frac{\tau_0}{\kappa},1\right]$ admits a unique fixed point which is  $a=1$, and then the fact that the same is true also in the interval 
$\left[1,1+\frac{\tau_0}{\kappa}\right]$ is a direct consequence of Remark \ref{amaggiore1}.

We first estimate $\lambda_a$ as in the proof of Proposition \ref{propositionstimal}: let $\varphi_a(x):=\frac{e^{-\frac{u_a(x)}{2}}}{\sqrt{\int_{-\pi}^{\pi} e^{-u_a(y)}dy}}=\sqrt{m_a(x)}$. Then, recalling that $1-a=\frac{\tau}{\kappa}$,  the couple $(\varphi_a, \lambda_a)$ solves
\begin{equation}\label{autova}
\begin{cases}
-\varphi_a''(x) + \frac{\tau  V(x)}2 \varphi_a(x) = \lambda_a \varphi_a(x) & \text{on $(-\pi, \pi)$}, \\
\varphi_a'(\pi) = \varphi_a'(-\pi) = 0\\ \int_{-\pi}^\pi \varphi_a^2(x)dx=1.
\end{cases}
\end{equation}
So, \[
\lambda_a = \inf_{\substack{\phi \in H^1(-\pi, \pi) \\ \int_{-\pi}^\pi \phi^2 = 1}} \int_{-\pi}^{\pi}  |\phi'|^2 + \frac{\tau V(x)}2 \phi^2 dx\leq  \inf_{\substack{\phi \in H^1(-\pi, \pi) \\ \int_{-\pi}^\pi \phi^2 = 1}} \int_{-\pi}^{\pi}  |\phi'|^2 + \frac{\tau x^2}4 \phi^2 dx 
\]
which yields \begin{equation}\label{elle} 0\leq \lambda_a\leq C\tau \end{equation}  as a straightforward consequence (just consider the constant competitor $\phi = \sqrt{1/2\pi}$). 

We multiply by $\phi_a$ the equation in \eqref{autova} and integrate by parts: we obtain, recalling \eqref{elle}, 
\begin{equation}\label{phi1}  \int_{-\pi}^\pi |\phi_a'|^2dx\leq \lambda_a \int_{-\pi}^\pi |\phi_a|^2dx=\lambda_a\leq C\tau.\end{equation} By the mean value theorem there exists $\xi\in [-\pi, \pi]$ such that $\frac{1}{2\pi}=\frac{1}{2\pi} \int_{-\pi}^\pi |\phi|^2 dx=  \phi^2(\xi)$. So we conclude, recalling that $m_a=\phi_a^2$, for all $x\in [-\pi, \pi]$ and using \eqref{phi1}
\begin{equation}\label{eqm1}  |m_a(x)-\frac1{2\pi}|= |m_a(x)-m_a(\xi)|=\left| \int_\xi^x 2\phi_a\phi_a'dx\right|\leq 2\left[\int_{-\pi}^\pi |\phi_a|^2dx\right]^{1/2}\left[\int_{-\pi}^\pi |\phi_a'|^2dx\right]^{1/2} \leq 2\sqrt{C\tau}. \end{equation} 
Again multiplying by $\phi_a''$ the equation in \eqref{autova} and integrating by parts, we get  by \eqref{elle}   \eqref{phi1} and the Young inequality that
\begin{eqnarray*} \int_{-\pi}^\pi |\phi_a''|^2dx&\leq&   \lambda_a \int_{-\pi}^\pi |\phi_a'|^2dx+\tau\int_{-\pi}^\pi V(x)\phi_a\phi_a''dx\\
&\leq &C^2\tau^2 +\frac{1}{2}\int_{-\pi}^\pi |\phi_a''|^2dx+\frac{\tau^2}{2}\int_{-\pi}^\pi  \|V\|_\infty^2|\phi_a|^2dx\leq (C^2+1)\tau^2 +\frac{1}{2}\int_{-\pi}^\pi |\phi_a''|^2dx .\end{eqnarray*}  
From this, recalling that $\phi'_a(\pm \pi)=0$ we conclude for all $x\in [-\pi, \pi]$, 
\[|\phi'_a(x)| =\left|\int_{-\pi}^x \phi_a''dx\right|\leq \sqrt{2\pi}\left(\int_{-\pi}^\pi |\phi_a''|^2dx\right)^{1/2}\leq C\tau. \]
Now we recall that $m_a(x)= e^{-u_a(x)}/ \int_{-\pi}^\pi e^{-u_a} dx$, and so for all $x\in [-\pi, \pi]$ and for $\tau>0$ sufficiently small such that $   m_a\geq 1/16$ (by \eqref{phi1}), we conclude 
\begin{equation}\label{stimaunuova} |u_a'(x)|=  \frac{|m_a'(x)|}{m_a(x)}= 2\frac{|\phi_a'(x)|}{\phi_a(x)}\leq 8C\tau,  \end{equation} and then  also $|u_a|\leq C\tau $.
 
By formula \eqref{derivata1} and by \eqref{phi1} we get that 
 \begin{multline}\label{stimatau1} |\lambda'_a+\kappa| =\left| -\kappa\int_{-\pi}^\pi V(x)m_a(x)dx +\kappa\right|=\left| -\kappa\int_{-\pi}^\pi V(x)\left(m_a(x)-\frac{1}{2\pi}\right)dx  \right| \\ \leq \kappa\int_{-\pi}^\pi V(x)\left|m_a(x)-\frac{1}{2\pi}\right|dx\leq \kappa C \sqrt{\tau}.  \end{multline}
 We consider now the function $v_a$ solution to \eqref{hjb1} with $v_a(0)=0$. Let us write $v_a(x)=\kappa (\cos x-1+z_a(x))$ for some function $z_a$.
 Then it is easy to check that $z_a$ is a solution to 
 \begin{equation}\label{hjb3} -z_a''+z_a' u_a' =\sin x \ u_a'-\frac{\kappa+\lambda_a'}{\kappa}
 \end{equation} with periodic boundary conditions and with $z_a(0)=0$. By the gradient estimates on $u_a$ \eqref{stimaunuova} and the estimate \eqref{stimatau1}, the right hand side of the previous equation is bounded by $C\sqrt{\tau}$ for some $C>0$, for $\tau<1$. It is a straightforward computation (by direct integration, and again by  the estimates on $u_a$) to show that $|z_a'(x)|\leq C\sqrt{\tau}$, and $|z_a(x)|\leq C\sqrt{\tau}$, for some $C>0$. 
 
Recalling  formula \eqref{derivata2} and the previous estimates on $z_a$, $m_a$  we compute $F_\kappa'(a)$: 
 \begin{eqnarray*} F'_\kappa(a)&=&  \frac{1}{\kappa} \int_{-\pi}^\pi (v_a'(y))^2m_a(y)dy= \kappa \int_{-\pi}^\pi ( z_a'(y)-\sin y)^2 m_a(y)dy\\ 
 &\geq &\kappa \left(\frac{1}{2\pi}-C\sqrt{\tau}\right) \int_{-\pi}^\pi ( z_a'(y)-\sin y)^2  dy \geq  \kappa \left(\frac{1}{2\pi}-C\sqrt{\tau}\right) (\pi -C\sqrt{\tau}).
\end{eqnarray*} In particular there exists $\tau_0=\tau_0(C)>0$ such that if  $\tau\leq \tau_0$, then  $F'_\kappa(a)=F'_\kappa\left(1-\frac{\tau}{\kappa}\right)\geq \frac{\kappa}{4}$.  
This implies immediately that  for all $\tau\in (0, \tau_0]$ there holds, for $\kappa>4$, 
\[F_\kappa\left(1-\frac{\tau}{\kappa}\right)\stackrel{  \tau'\in (0,\tau)}{=} F_\kappa(1)-\frac{\tau}{\kappa}F'_\kappa\left(1-\frac{\tau'}{\kappa}\right)
\stackrel{\tau'<\tau_0}{\leq} 1-\frac{\tau}{\kappa}\frac{\kappa}{4}=1-\frac{\tau}{4}<1-\frac{\tau}{\kappa}.\]
\end{proof} 

Using the previous results, we conclude with the proof of the main result of this section.
 
 \begin{proof}[Proof of Theorem \ref{ex2}]  
 
  Let $\delta=\frac{\tau_0}{4}$, where $\tau_0$ is as in Proposition \ref{propositiontau}. By Theorem \ref{ex1}, there exists $\kappa_0:=\kappa_0(\tau_0)>4$ (using the same notation as in Theorem \ref{ex1}) such that
 for   $\kappa\geq \kappa_0$  the map $F_\kappa:\left[0, 1-\frac{\tau_0}{4}\right]\to \left[0, 1-\frac{\tau_0}{4}\right]$ admits a unique fixed point $\bar a$. 
 
 By Proposition  \ref{propositiontau}   for all $\kappa\geq \kappa_0$ there holds $F_\kappa\left(1-\frac{\tau_0}{\kappa}\right)\leq   1-\frac{\tau_0}{4}$. Since $F_\kappa$ is a nondecreasing map by \eqref{derivata2}, this implies that $F_\kappa(a)<a$ for all $a\in \left[1-\frac{\tau_0}{4}, 1-\frac{\tau_0}{\kappa}\right]$. 
 
Finally, again by Proposition \ref{propositiontau},  $F_\kappa$ in $\left[1-\frac{\tau_0}{\kappa},1\right]$ admits a unique fixed point which is  $a=1$. 
This implies that there exists a unique fixed point $\bar a\in [0,1)$.  Note that by Remark \ref{amaggiore1}, $2-\bar a$ is the unique fixed point in $(1,2]$. 
  \end{proof}

\section{Local dynamic stability of the self-organizing solution}\label{sezdinamica}

Let $(\bar u, \bar \lambda, \bar m)$ be the unique stationary even self-organizing solution with $\int_{-\pi}^\pi \cos x\  \bar m(x)dx>0$, which has been obtained in the previous section, under the assumption that $\kappa \ge \kappa_0$ (see Theorem \ref{ex2}). We show in this section a \textit{local} exponential stability property of $(\bar u, \bar \lambda, \bar m)$. 

We consider the  dynamic solutions $( \tilde u ,m)$ of \eqref{mfgev1}. First of all we observe that if we define 
\begin{equation}\label{tra} u(x,t):=\tilde u(x,t)-\kappa  \int_0^t \int_{-\pi}^{\pi}  \cos y\  m(s,y)dyds\end{equation} then $(u,m)$ is a dynamic solution to   \begin{equation}\label{mfgev}  \begin{cases} -u_t-u_{xx}+\frac{1}{2}|u_x|^2=   \kappa V(x) \left[1-\int_{-\pi}^{\pi} V(y)m(t,y)dy\right]  
\\ m_t-m_{xx}-(mu_x)_x=0\\ m_x(t, \pi )=m_x(t, -\pi )=0\qquad u_x(t, \pi )=u_x(t, -\pi)=0\\ 
\text{$u(\cdot,t),m(\cdot,t)$ are even, $\int_{-\pi}^{\pi} m(x,t)dx=1, m(\cdot,t) \ge 0$ for all $t$}.\end{cases}   \end{equation}
Note that Neumann boundary conditions at the boundary of $[-\pi, \pi]$ are equivalent to requiring $(u,m)$ to be $2\pi$-periodic.

We are going to show that if $(u,m)$ is a  solution to \eqref{mfgev} such that the    density and the  optimal control $(m, u_x)$ remain  for all the time in a suitable neighborhood of the equilibrium density and of the ergodic optimal control $(\bar m,\bar u_x)$, then $(m,u_x)$ is going to   converge exponentially fast to $(\bar m,\bar u_x)$ as $T$ goes to infinity. In particular, this will also imply that the associated solution $(\tilde u, m)$ to \eqref{mfgev}, according to \eqref{tra}, satisfies the same exponential stability property. 

  We introduce the following constant 
\begin{equation}\label{defQ1} Q:= Q_\kappa=\int_{-\pi}^\pi  \kappa x^4 \bar m(x)dx.\end{equation} It is easy to check, by using the upper and lower bounds \eqref{mbound} obtained in Corollary \ref{corollary2.9},  that for $\kappa>\kappa_0$, $Q$ can be controlled above and below by some positive constants independent of $\kappa$ (in fact, depending on $c_1$, $c_2$, $c_3$ in \eqref{mbound}).  

Our main result reads as follows.

\begin{theorem} \label{conv1} Let $(u,m)$ be a solution of \eqref{mfgev}. Assume that
 \begin{equation}\label{initialdataun} 0<m(t,x)\leq c \kappa^{\frac{1}{4}}  \bar m(x) 
\end{equation}  for all $x\in [-\pi,\pi]$, $t\in [0, T]$, where $c$ is specified below (see \eqref{choicec}).
Then there exist $\kappa_1\geq \kappa_0 \wedge 1$, 
and  $C>0$ (independent of $\kappa, T$), such that for all $0 \le t \leq T $ and $\kappa\geq \kappa_1$, there holds 
\begin{equation}\label{convave}
\dashint_{t}^{t+C\kappa^{-1/4}}  \int_{-\pi}^{\pi } |u_x(\tau, x)- \bar u_x(x)|^2 \bar m   +  Q\frac{|m(\tau, x)-\bar m(x)|^2}{ \bar m(x)} dx d\tau \le K (e^{-\omega t} + e^{-\omega (T-t)}),
\end{equation}
where
\[
K = 4\sum_{t=0, t=T} \int_{-\pi}^{\pi } |u_x(x,t)- \bar u_x(x)|^2 \bar m(x)  +  Q\frac{|m(x, t)-\bar m(x)|^2}{ \bar m(x)} dx , \qquad \omega = \frac{\log 2}{C} \kappa^{1/4}.
\]
Moreover, for every $t \in [0,T]$,
\begin{equation}\label{convm}
 \int_{-\pi}^{\pi }  \frac{|m(t, x)-\bar m(x)|^2}{ \bar m(x)} dx \le K C' (\kappa^{1/4} + 1) (e^{-\omega t} + e^{-\omega (T-t)}).
\end{equation}
\end{theorem}

The constant $C$ will actually depend on $C_P$ (see \eqref{poincare}  below) and $Q$.

\begin{remark}\normalfont The way the previous statement quantifies the convergence of  $u_x$ to $\bar u_x$ is in time average (note that the length of time integration vanishes as $\kappa$ goes to infinity). One could get an exponential convergence \textit{pointwise} in time, as it is done for $m$ to $\bar m$ in \eqref{convm}, by coupling \eqref{convave} with suitable estimates on the linearized HJ equation, which are not developed here. 

Another way would be to use Lemma \ref{poincinterp} which yields pointwise information in time right away from \eqref{convave}, though it involves a (uniform in $T$) Lipschitz control on $t \mapsto \int |u_x(t, x)- \bar u_x(x)|^2 \bar m   +  \frac{|m(t, x)-\bar m(x)|^2}{ \bar m(x)} dx$. Such control can be derived, but it should be quite sensitive to the value of $\kappa$.
\end{remark}

To prove Theorem \ref{conv1}, we first rescale the problem as in Remark \ref{rescaled}; let us consider 
\[w(t,x)= u(t\kappa^{-\frac{1}{2}}, x\kappa^{-\frac14}), \qquad \mu(t,x)= \kappa^{-\frac14} m(t\kappa^{-\frac{1}{2}}, x\kappa^{-\frac14}),\]
that satisfy in the rescaled space-time cylinder $(0,T\kappa^{\frac{1}{2}})\times [-\pi\kappa^{\frac14}, \pi\kappa^{\frac14}]$ the system
  \begin{equation}\label{mfgresc}  \begin{cases} -w_t-w_{xx}+\frac{1}{2}|w_x|^2=  V_\kappa(x) 
  \left[1-\kappa^{-\frac{1}{2}}\int_{-\pi\kappa^{\frac14}}^{\pi\kappa^{\frac14}} V_\kappa(y)\mu(y)dy\right]
\\ \mu_t-\mu_{xx}-(\mu w_x)_x=0\\ \mu_x(t, \pi\kappa^{\frac14})=\mu_x(t, -\pi\kappa^{\frac14})=0\qquad w_x(t, \pi\kappa^{\frac14})=w_x(t, -\pi\kappa^{\frac14})=0\\ 
\text{$w(\cdot,t),\mu(\cdot,t)$ are even, $\int_{-\pi \kappa^{\frac14}}^{\pi \kappa^{\frac14}} \mu(x,t)dx=1, \mu(\cdot,t) \ge 0$ for all $t$.}
\end{cases}
 \end{equation}

   Let  $\kappa>\kappa_0$  and $(\bar w, \bar \mu,\bar \lambda)$ be the rescaling, according to Remark \ref{rescaled}, of the unique even self organizing solution with $ \bar \lambda>0 $ obtained in Theorem \ref{ex2}:
   \[\bar w(x)= \bar u(  x\kappa^{-\frac14}), \qquad \bar\mu(x)= \kappa^{-\frac14} \bar m(  x\kappa^{-\frac14}).\]   We recall that by Proposition \ref{uabobelo} and Corollary \ref{corollary2.9}, we get that for some constants independent of $\kappa$,
  \begin{equation}\label{barubarm} c_1|x|^2-c_3\leq \bar w(x)\leq c_2|x|^2+c_3  \qquad c_4e^{-c_2 x^2}\leq \bar \mu(x)\leq c_5e^{-c_1 x^2 } \qquad x\in[-\kappa^{\frac{1}{4}}\pi, \kappa^{\frac14}\pi]. \end{equation} 
Recall also that we fixed $\bar w(0)=0$ and moreover there holds $\bar w(0)=\min \bar w$, by simmetry of $w$. 
Due to \eqref{barubarm}, the following  weighted Poincar\'e  inequality holds.
\begin{theorem}[Poincar\'e weighted inequality] \label{poincarethm}  Let $(\bar w, \bar \mu,\bar \lambda)$ as in \eqref{barubarm} and $\kappa\geq \kappa_0$, as in 
Theorem \ref{ex2}. Then there exist $\kappa_1\geq \kappa_0$ and  a constant $C_P$ independent of $\kappa$ such that for all  $\kappa\geq \kappa_1$, $f\in H^1_{\bar\mu} (-\kappa^{\frac14}\pi, \kappa^{\frac14}\pi)$ with 
$\int_{-\kappa^{\frac14}\pi}^{\kappa^{\frac14}\pi} f(x)\bar \mu(x)dx=0$ there holds
\begin{equation} \label{poincare} \int_{-\kappa^{\frac14}\pi}^{\kappa^{\frac14}\pi} f^2(x)\bar \mu(x)dx\leq C_P \int_{-\kappa^{\frac14}\pi}^{\kappa^{\frac14}\pi} 
f_x^2(x)\bar \mu(x)dx.\end{equation}
\end{theorem} 
The proof is reported for completeness in the appendix. Note that the constant $Q$ introduced in \eqref{defQ1} coincides with  \begin{equation}\label{defQ} Q=  \int_{-\pi\kappa^{\frac14}}^{\pi\kappa^{\frac14}}x^4\bar\mu(x)dx.\end{equation}

 We can now specify $c$ in the statement of Theorem \ref{conv1}
 \begin{equation}\label{choicec}
c = \sqrt{\frac{C_P}{Q}}
 \end{equation}
 so that \eqref{initialdataun} reads
  \begin{equation}\label{initialdata} 0<\mu(t,x)\leq \kappa^{\frac{1}{4}} \sqrt{\frac{C_P}{Q}} \bar \mu(x) 
\end{equation}  for all $x\in [-\pi\kappa^{\frac14}, \pi\kappa^{\frac14}]$, $t\in [0, T\kappa^{\frac12}]$.



 Let us define  $\zeta(t,x)= \mu(t,x)-\bar \mu(x)$ and $v(t,x)=w(t,x)-\bar w(x)-\bar \lambda (T-t)$. They are solutions to 
   \begin{equation}\label{lin}  \begin{cases} -v_t-v_{xx}+\frac{1}{2}|v_x|^2+ \bar w_x v_x= 
   - \kappa^{-\frac{1}{2}}V_\kappa(x)\int_{-\pi\kappa^{\frac14}}^{\pi\kappa^{\frac14} } V_\kappa(y)\zeta(t,y)dy 
\\ \zeta_t-\zeta_{xx}-(\zeta \bar w_x)_x=(\mu v_x)_x  \\ \zeta(0,x)=\kappa^{-\frac14} m_0(x\kappa^{-\frac14})( x)-\bar \mu(x),\ \  v(T\kappa^{\frac12},x)=  u_T(x\kappa^{-\frac14})-\bar w(x)\end{cases}   \end{equation}    
with Neumann boundary conditions. Observe that $\int_{-\pi\kappa^{\frac14}}^{\pi\kappa^{\frac14}} \zeta(t,x)dx=0$ for all $t$.

 The main result that will provide the exponential convergence is the following.  \begin{proposition} \label{expo} 
Let 
\[\Phi(t):=  \int_{-\pi\kappa^{\frac14}}^{\pi\kappa^{\frac14} } \bar \mu(x) |v_x(t,x)|^2 dx  + \frac{Q}{ \kappa^{1/2}} \int_{-\pi\kappa^{\frac14}}^{\pi\kappa^{\frac14} } \frac{|\zeta(t,x)|^2}{\bar \mu(x)} dx.\]
Let us assume that \eqref{initialdata} holds. Then there exists $\kappa_1\geq \kappa_0 $, where $\kappa_0$ is  as in Theorem \ref{ex2},  such that  for all $0\leq t_1<t_2\leq T \kappa^{1/2}$ and $\kappa\geq \kappa_1$, there holds 
\begin{equation}\label{claim} \int_{t_1}^{t_2} \Phi(t)dt \leq 4 \left(C_P+\frac{1}{C_P}+\frac{1}{Q}\right) \kappa^{1/4} \big(\Phi(t_1)+\Phi(t_2)\big).
\end{equation}
\end{proposition} 

 For the proof of this Proposition we need some lemmata.   
First of all we have the following result, obtained by duality arguments. 
\begin{lemma}\label{lemma1} \begin{align}\label{dualitymuv1} \int_{t_1}^{t_2}  \int_{-\pi\kappa^{\frac14}}^{\pi\kappa^{\frac14} } \bar \mu(x) v_x^2  dxdt & -2\kappa^{-\frac12}\int_{t_1}^{t_2} \left(\int_{-\pi\kappa^{\frac14}}^{\pi\kappa^{\frac14} }V_k(x)\zeta(t,x)dx\right)^2\\ \nonumber & \leq  C_P \kappa^{1/4} \int_{-\pi\kappa^{\frac14}}^{\pi\kappa^{\frac14} } \bar \mu(x)v_x^2(t_1,x)dx
 + \frac1{\kappa^{1/4}} \int_{-\pi\kappa^{\frac14}}^{\pi\kappa^{\frac14} } \frac{|\zeta(t_1,x)|^2}{\bar \mu(x)}dx\\  \nonumber &+  C_P \kappa^{1/4} \int_{-\pi\kappa^{\frac14}}^{\pi\kappa^{\frac14} } \bar \mu(x)v_x^2(t_2,x)dx
 +  \frac1{\kappa^{1/4}} \int_{-\pi\kappa^{\frac14}}^{\pi\kappa^{\frac14} } \frac{|\zeta (t_2,x)|^2}{\bar \mu(x)}dx. \end{align} 
\end{lemma} \begin{proof}
By duality, since $v,\zeta$ are solutions to \eqref{lin}, and recalling that $\mu>0$,  we get
\begin{eqnarray*} &&\frac{d}{dt}\int_{-\pi\kappa^{\frac14}}^{\pi\kappa^{\frac14} } v(t,x)\zeta(t,x) dx=-\int_{-\pi\kappa^{\frac14}}^{\pi\kappa^{\frac14} } \frac{\mu(t,x)+\bar \mu(x)}{2} v_x^2 (t,x) dx+\kappa^{-\frac12}\left(\int_{-\pi\kappa^{\frac14}}^{\pi\kappa^{\frac14} }V_k(x)\zeta (t,x)dx\right)^2\\ &\leq& -\int_{-\pi\kappa^{\frac14}}^{\pi\kappa^{\frac14} } \frac{\bar \mu(x)}{2} v_x^2(t,x)  dx+\kappa^{-\frac12}\left(\int_{-\pi\kappa^{\frac14}}^{\pi\kappa^{\frac14} }V_k(x)\zeta(t,x)dx\right)^2. \end{eqnarray*}
 Integrating in $(t_1, t_2)\subseteq(0,T \kappa^{1/2})$ we get
\begin{eqnarray}\label{dualitymuv}&& \int_{t_1}^{t_2} \int_{-\pi\kappa^{\frac14}}^{\pi\kappa^{\frac14} } \frac{\bar \mu(x)}{2} v_x^2(t,x)  dxdt-\kappa^{-\frac12}\int_{t_1}^{t_2} \left(\int_{-\pi\kappa^{\frac14}}^{\pi\kappa^{\frac14} }V_k(x)\zeta(t,x)dx\right)^2\\ \nonumber & \leq & \int_{-\pi\kappa^{\frac14}}^{\pi\kappa^{\frac14} } \zeta(t_1,x)  v(t_1,x) dx -\int_{-\pi\kappa^{\frac14}}^{\pi\kappa^{\frac14} } \zeta(t_2, x) v(t_2,x) dx .  \nonumber \end{eqnarray} 
 
Let $t\in [0,T\kappa^{1/2}]$ and $c(t)= \int_{-\pi\kappa^{\frac14}}^{\pi\kappa^{\frac14} }   v(t,x) \bar\mu(x) dx$. Since $\int_{-\pi\kappa^{\frac14}}^{\pi\kappa^{\frac14} } \zeta(t, x)  dx=0$ for all $t$, there holds, 
also using Young inequality and the Poincar\'e inequality \eqref{poincare}:
\begin{eqnarray*} \left|\int_{-\pi\kappa^{\frac14}}^{\pi\kappa^{\frac14} } \zeta(t, x) v(t,x)  dx\right|&=&\left| \int_{-\pi\kappa^{\frac14}}^{\pi\kappa^{\frac14} } \zeta (t,x)( v(t,x)-c(t))dx\right| \\
&\leq &\frac{\kappa^{1/4}}2 \int_{-\pi\kappa^{\frac14}}^{\pi\kappa^{\frac14} } \bar \mu(x)( v(t,x)-c(t))^2 dx +\frac{1}{2 \kappa^{1/4} } \int_{-\pi\kappa^{\frac14}}^{\pi\kappa^{\frac14} } \frac{|\zeta(t,x)|^2}{\bar \mu(x)}dx\\
&\leq& \frac{C_P \kappa^{1/4} }{2}   \int_{-\pi\kappa^{\frac14}}^{\pi\kappa^{\frac14} } \bar \mu(x)v_x^2(t,x)dx
 +\frac{1}{2\kappa^{1/4}} \int_{-\pi\kappa^{\frac14}}^{\pi\kappa^{\frac14} } \frac{|\zeta(t,x)|^2}{\bar \mu(x)}dx.\end{eqnarray*}
 Substituting this estimate in \eqref{dualitymuv} per $t=t_1,t_2$  we conclude that \eqref{dualitymuv1} holds. \end{proof}

\begin{lemma}\label{lemmavkappa} For all $0\leq t_1<t_2\leq T\kappa^{\frac12}$, there holds 
\begin{equation}\label{vkappa}\int_{t_1}^{t_2} \left(\int_{-\pi\kappa^{\frac14}}^{\pi\kappa^{\frac14} }V_k(x)\zeta (t,x)dx\right)^2dt\leq \frac{ Q }{4} \int_{t_1}^{t_2} \int_{-\pi\kappa^{\frac14}}^{\pi\kappa^{\frac14} }  \frac{|\zeta(t,x)|^2}{\bar \mu(x)} dx.\end{equation} \end{lemma}
\begin{proof} 
First of all we observe, recalling \eqref{Vkbound}, that \[\int_{t_1}^{t_2}\left(\int_{-\pi\kappa^{\frac14}}^{\pi\kappa^{\frac14} }V_k(x)\zeta(t,x)dx\right)^2dt \leq 
\frac14 \int_{t_1}^{t_2}\left(\int_{-\pi\kappa^{\frac14}}^{\pi\kappa^{\frac14} } |\zeta(t,x)|x^2 dx\right)^2dt\]
for some constant $C$ not depending on $\kappa$.   By H\"older inequality and \eqref{barubarm} we get
\[\left(\int_{-\pi\kappa^{\frac14}}^{\pi\kappa^{\frac14} } |\zeta(t,x)|x^2 dx\right)^2\leq  
\left(\int_{-\pi\kappa^{\frac14}}^{\pi\kappa^{\frac14} } \frac{|\zeta(t,x)|^2}{\bar \mu(x)} dx\right)
  \left(\int_{-\pi\kappa^{\frac14}}^{\pi\kappa^{\frac14} }  x^4 \bar \mu(x) dx\right) .\] Substituting in the previous inequality we get the conclusion.
\end{proof}
 \begin{lemma}\label{lemmabarm} Let assume that \eqref{initialdata} holds.   For all  $0\leq t_1<t_2\leq T\kappa^{1/2}$, there holds 
\begin{multline} \int_{t_1}^{t_2} \int_{-\pi\kappa^{\frac14}}^{\pi\kappa^{\frac14} } \frac{|\zeta(t,x)|^2}{\bar \mu(x)}  dxdt  +\frac{1}{C_P} \int_{-\pi\kappa^{\frac14}}^{\pi\kappa^{\frac14} }\frac{|\zeta(t_2,x)|^2}{\bar \mu(x)}dx\\ 
\leq \frac{ \kappa^{\frac{1}{2}}}{Q}\int_{t_1}^{t_2}  \int_{-\pi\kappa^{\frac14}}^{\pi\kappa^{\frac14} }v_x^2(t,x)\bar \mu(x) dxdt +\frac{1}{C_P} \int_{-\pi\kappa^{\frac14}}^{\pi\kappa^{\frac14} }\frac{|\zeta(t_1,x)|^2}{\bar \mu(x)}dx. \label{muquadro} 
\end{multline} 
\end{lemma}
\begin{proof} Note that the equation satisfied by $\mu$ in \eqref{mfgresc} can be written as
\[\mu_t-\mu_{xx}-(\mu\bar w_x)_x= (\mu v_x)_x,\] multiply it by $\frac{\mu(t,x)}{\bar \mu(x)}-1=\frac{\zeta(t,x)}{\bar \mu(x)}$ and integrate in $[-\pi\kappa^{-\frac14}, \pi\kappa^{\frac14}]$
\begin{eqnarray*}&&
\frac{1}{2}\frac{d}{dt}\int_{-\pi\kappa^{\frac14}}^{\pi\kappa^{\frac14} } \bar \mu(x) \left(\frac{\mu(t,x)}{\bar \mu(x)}-1\right)^2dx+\int_{-\pi\kappa^{\frac14}}^{\pi\kappa^{\frac14} } (\mu_x(t,x)+\bar w_x(x)\mu(t,x)) \left(\frac{\mu(t,x)}{\bar \mu(x)}-1\right)_x dx\\ &=&-  \int_{-\pi\kappa^{\frac14}}^{\pi\kappa^{\frac14} }v_x(t,x) \mu(t,x) \left(\frac{\mu(t,x)}{\bar \mu(x)}-1\right)_xdx.\end{eqnarray*}
 Recalling that $\bar \mu(x)= e^{-\bar w(x)+c}$, see \eqref{hjbresc}, we get $\bar \mu_x(x)=-\bar w_x(x) \bar \mu(x)$ and so 
 \[ \left(\frac{\mu(t,x)}{\bar \mu(x)}-1\right)_x= \frac{\mu_x(t,x)+\bar w_x(x) \mu(t,x)}{\bar \mu(x)}.\]
 Substituting in the previous equality and using the Young inequality for the right hand side we obtain 
 \begin{eqnarray*}&& 
\frac{1}{2}\frac{d}{dt}\int_{-\pi\kappa^{\frac14}}^{\pi\kappa^{\frac14} } \bar \mu(x) \left(\frac{\mu(t,x)}{\bar \mu(x)}-1\right)^2dx
+\int_{-\pi\kappa^{\frac14}}^{\pi\kappa^{\frac14} } \bar \mu(x) \left[\left(\frac{\mu(t,x)}{\bar \mu(x)}-1\right)_x\right]^2dx\\
&=&-  \int_{-\pi\kappa^{\frac14}}^{\pi\kappa^{\frac14} }v_x(t,x) \mu(t,x) \left(\frac{\mu(t,x)}{\bar \mu(x)}-1\right)_xdx\\
&\leq& \frac{1}{2} \int_{-\pi\kappa^{\frac14}}^{\pi\kappa^{\frac14} }v_x^2(t,x)\frac{\mu^2(t,x)}{\bar \mu(x)}dx+\frac{1}{2}  \int_{-\pi\kappa^{\frac14}}^{\pi\kappa^{\frac14} } \bar \mu(x)\left[\left(\frac{\mu(t,x)}{\bar \mu(x)}-1\right)_x\right]^2dx. \end{eqnarray*}
Therefore we get, recalling that by \eqref{initialdata} $\mu^2(t,x)\leq  \kappa^{\frac12} \frac{C_P}{Q} \bar\mu^2(x)$, 
 \begin{eqnarray*}&& 
 \frac{d}{dt}\int_{-\pi\kappa^{\frac14}}^{\pi\kappa^{\frac14} } \bar \mu(x) \left(\frac{\mu(t,x)}{\bar \mu(x)}-1\right)^2dx
+\int_{-\pi\kappa^{\frac14}}^{\pi\kappa^{\frac14} } \bar \mu(x) \left(\frac{\mu(t,x)}{\bar \mu(x)}-1\right)^2_x dx 
\leq   \kappa^{\frac12} \frac{C_P}{Q}  \int_{-\pi\kappa^{\frac14}}^{\pi\kappa^{\frac14} }v_x^2(t,x) \bar \mu(x)dx.  \end{eqnarray*}
By the Poincar\'e inequality, \eqref{poincare}, since  $\int_{-\pi\kappa^{\frac14}}^{\pi\kappa^{\frac14} } \bar \mu(x) \left(\frac{\mu(t,x)}{\bar \mu(x)}-1\right)dx=0$, 
we get
 \begin{eqnarray*} 
 \frac{d}{dt}\int_{-\pi\kappa^{\frac14}}^{\pi\kappa^{\frac14} } \bar \mu(x) \left(\frac{\mu(t,x)}{\bar \mu(x)}-1\right)^2dx +C_P\int_{-\pi\kappa^{\frac14}}^{\pi\kappa^{\frac14} } \bar \mu(x) \left(\frac{\mu(t,x)}{\bar \mu(x)}-1\right)^2dx
\leq  \kappa^{\frac12} \frac{C_P}{Q} \int_{-\pi\kappa^{\frac14}}^{\pi\kappa^{\frac14} }v_x^2(t,x) \bar \mu(x) dx.  \end{eqnarray*}
By integration on $(t_1, t_2)$ we get, recalling that  $\bar \mu(x) \left(\frac{\mu(t,x)}{\bar \mu(x)}-1\right)^2=\frac{|\zeta(t,x)|^2}{\bar \mu(x)}$, 
 \begin{multline*} \frac{1}{C_P}  \int_{-\pi\kappa^{\frac14}}^{\pi\kappa^{\frac14} }\frac{|\zeta(t_2,x)|^2}{\bar \mu(x)}dx+  \int_{t_1}^{t_2} \int_{-\pi\kappa^{\frac14}}^{\pi\kappa^{\frac14} }\frac{|\zeta(t,x)|^2}{\bar \mu(x)} dx\\ \leq \frac{\kappa^{\frac12}}{Q} \int_{t_1}^{t_2}   \int_{-\pi\kappa^{\frac14}}^{\pi\kappa^{\frac14} }v_x^2(t,x) \bar \mu(x) dx +\frac{1}{C_P}  \int_{-\pi\kappa^{\frac14}}^{\pi\kappa^{\frac14} }\frac{|\zeta(t_1,x)|^2}{\bar \mu(x)}dx  . \end{multline*} 
 \end{proof} 
 We are ready to prove Proposition \ref{expo}.
 \begin{proof}[Proof of Proposition \ref{expo}]
 We rewrite the inequality   \eqref{dualitymuv1} by recalling the definition of $\Phi(t)$: 
\begin{equation}\label{dual}  \int_{t_1}^{t_2}  \int_{-\pi\kappa^{\frac14}}^{\pi\kappa^{\frac14} } \bar \mu(x) v_x^2(t,x)  dxdt-2\kappa^{-\frac12}\int_{t_1}^{t_2} \left(\int_{-\pi\kappa^{\frac14}}^{\pi\kappa^{\frac14} }V_k(x)\zeta(t,x)dx\right)^2\leq \kappa^{1/4} \left(C_P+\frac1Q\right)\big(\Phi(t_1)+\Phi(t_2)\big).\end{equation} 
By \eqref{vkappa} and \eqref{muquadro} we obtain 
\begin{multline*} -2\kappa^{-\frac12}\int_{t_1}^{t_2} \left(\int_{-\pi\kappa^{\frac14}}^{\pi\kappa^{\frac14} }V_k(x)\zeta(t,x)dx\right)^2\geq  
 -\kappa^{-\frac12}\frac{Q}{2}  
 \int_{t_1}^{t_2} \int_{-\pi\kappa^{\frac14}}^{\pi\kappa^{\frac14} }  \frac{|\zeta(t,x)|^2}{\bar \mu(x)} dx\\ 
 \geq - \frac{1}{2}   \int_{t_1}^{t_2}  \int_{-\pi\kappa^{\frac14}}^{\pi\kappa^{\frac14} }v_x^2(t,x) \bar \mu(x) dxdt 
  -  \kappa^{-\frac12}\frac{Q}{2C_P} \int_{-\pi\kappa^{\frac14}}^{\pi\kappa^{\frac14} }\frac{|\zeta(t_1,x)|^2}{\bar \mu(x)}dx.   \end{multline*} 
Now, using \eqref{dual}, the previous  inequality and again \eqref{muquadro} we get 
\begin{multline*} \kappa^{1/4 } (C_P+1/Q)\big(\Phi(t_1)+\Phi(t_2)\big)\geq \\ \geq  \int_{t_1}^{t_2}  \int_{-\pi\kappa^{\frac14}}^{\pi\kappa^{\frac14} } \bar \mu(x) v_x^2(t,x)  dxdt-2\kappa^{-\frac12}\int_{t_1}^{t_2} \left(\int_{-\pi\kappa^{\frac14}}^{\pi\kappa^{\frac14} }V_k(x)\zeta(t,x)dx\right)^2\\ \geq\frac12   \int_{t_1}^{t_2}  \int_{-\pi\kappa^{\frac14}}^{\pi\kappa^{\frac14} } \bar \mu(x) v_x^2(t,x)  dxdt-\kappa^{-\frac12}\frac{Q}{2C_P} \int_{-\pi\kappa^{\frac14}}^{\pi\kappa^{\frac14} }\frac{|\zeta(t_1,x)|^2}{\bar \mu(x)}dx  \\
\geq \frac14   \int_{t_1}^{t_2}  \int_{-\pi\kappa^{\frac14}}^{\pi\kappa^{\frac14} } \bar \mu(x) v_x^2(t,x)  dxdt + \frac{Q}{4\kappa^{1/2}}\int_{t_1}^{t_2} \int_{-\pi\kappa^{\frac14}}^{\pi\kappa^{\frac14} } \frac{|\zeta(t,x)|^2}{\bar \mu(x)} dx - \\
-\frac1{\kappa^{1/2}}   \left(\frac{Q}{4C_P} +\frac{Q}{2C_P}\right)  \int_{-\pi\kappa^{\frac14}}^{\pi\kappa^{\frac14} }\frac{|\zeta(t_1,x)|^2}{\bar \mu(x)}dx  \end{multline*}
 and so, by the definition of $\Phi$, 
\[
\kappa^{1/4}(C_P+1/Q)(\Phi(t_1)+\Phi(t_2))\geq \frac14 \int_{t_1}^{t_2} \Phi(t) dt - \frac1{C_P} \Phi(t_1).
 \]
   \end{proof}
   We conclude with the proof of Theorem \ref{conv1}.
 
    \begin{proof}[Proof of Theorem \ref{conv1}] Let $C=16 \left(C_P+\frac{1}{C_P}+\frac1Q\right)$. A combination of Proposition \ref{expo} and Lemma \ref{expdecaylemma} yields
\begin{equation}\label{fi}
\frac1{C \kappa^{1/4}} \int_{t}^{t+ C \kappa^{1/4}} \Phi(s) ds \le 4  (e^{-\tilde \omega t} + e^{-\tilde\omega (T\kappa^{1/2}-t)})[\Phi(0) + \Phi(T\kappa^{\frac12})], \qquad \text{where \, $\tilde\omega = \frac{\log 2}{C \kappa^{1/4}}$}.
\end{equation}
Since
\[
\Phi(s) = \frac1{\kappa^{1/2}} \int_{-\pi}^{\pi } |u_x(x,s \kappa^{-1/2})- \bar u_x(x)|^2 \bar m(x)  +  Q\frac{|m(x, s \kappa^{-1/2})-\bar m(x)|^2}{ \bar m(x)} dx,
\]
we get in \eqref{fi} by performing a change of variables $\tau = s \kappa^{-1/2}$ 
\begin{eqnarray*} && \frac1{C \kappa^{1/4}}
 \int_{t}^{t+ C \kappa^{1/4}} \frac1{\kappa^{1/2}} \int_{-\pi}^{\pi } |u_x(x,s \kappa^{-1/2})- \bar u_x(x)|^2 \bar m(x)  +  Q\frac{|m(x, s \kappa^{-1/2})-\bar m(x)|^2}{ \bar m(x)} dx ds \\
 &=&  \frac1{C \kappa^{1/4}}
 \int_{t\kappa^{-1/2}}^{t\kappa^{-1/2}+ C \kappa^{-1/4}}  \int_{-\pi}^{\pi } |u_x(x,\tau)- \bar u_x(x)|^2 \bar m(x)  +  Q\frac{|m(x,\tau)-\bar m(x)|^2}{ \bar m(x)} dx d\tau\\
 &\leq&    (e^{-\tilde\omega\kappa^{1/2} t\kappa^{-1/2}} + e^{-\tilde\omega\kappa^{1/2} (T-t\kappa^{-1/2})})\frac{K}{\kappa^{1/2}}.
\end{eqnarray*} 
Replacing now $t\kappa^{-1/2}$ by $t$ and observing that  $\tilde \omega\kappa^{1/2}=\omega$,  we obtain the first assertion.

Applying now the Mean Value Theorem in \eqref{fi}, for every $t \in [C \kappa^{1/4}, T \kappa^{1/2}]$ there exists $\xi =\xi(t) \in [t-C \kappa^{1/4}, t]$ such that
\[
\frac{Q}{ \kappa^{1/2}} \int_{-\pi\kappa^{\frac14}}^{\pi\kappa^{\frac14} } \frac{|\zeta|^2}{\bar \mu}(\xi, x) dx \le \Phi(\xi) \le 4  (e^{-\tilde \omega (t-C\kappa^{1/4})} + e^{-\tilde\omega (T\kappa^{1/2}-t+C\kappa^{1/4})})[\Phi(0) + \Phi(T\kappa^{\frac12})].
\]
By Lemma \ref{lemmabarm},
\begin{align*}
\int_{-\pi\kappa^{\frac14}}^{\pi\kappa^{\frac14} }\frac{|\zeta|^2}{\bar \mu}(t,x)dx
&\leq \frac{C_P  \kappa^{\frac{1}{2}}}{Q}\int_{\xi}^{t}  \int_{-\pi\kappa^{\frac14}}^{\pi\kappa^{\frac14} }v_x^2\bar \mu dxdt + \int_{-\pi\kappa^{\frac14}}^{\pi\kappa^{\frac14} }\frac{|\zeta|^2}{\bar \mu}(\xi,x)dx\\
&\le \frac{C_P  \kappa^{\frac{1}{2}}}{Q}\int_{t-C\kappa^{1/4}}^{t} \Phi(t) dt + \int_{-\pi\kappa^{\frac14}}^{\pi\kappa^{\frac14} }\frac{|\zeta|^2}{\bar \mu}(\xi,x)dx, \label{muquadro} 
\end{align*}
hence using again \eqref{fi} we get
\[
\int_{-\pi\kappa^{\frac14}}^{\pi\kappa^{\frac14} }\frac{|\zeta|^2}{\bar \mu}(t,x)dx \le 4 \frac{\kappa^{1/2}}Q \left(C\kappa^{1/4}+ 1 \right)  (e^{-\tilde \omega (t-C\kappa^{1/4})} + e^{-\tilde\omega (T\kappa^{1/2}-t+C\kappa^{1/4})})[\Phi(0) + \Phi(T\kappa^{\frac12})].
\]
Going now back to the original space-time scale we obtain the statement.
\end{proof}   

   \appendix
   \section{Some estimates and a Poincar\'e weighted inequality}
   
    \begin{lemma}\label{expdecaylemma} Assume that $\Phi : [0,T] \to [0, \infty)$ satisfies
 \[
 \int_{t_1}^{t_2} \Phi(s) d s \le C [\Phi(t_1) + \Phi(t_2)]
 \]
 for some $C > 0$ and all $t_1 < t_2 \in [0,T]$ (assume also that $T \ge 8C$). Then, for all $t \in [0, T-4C]$,
 \[
 \frac{1}{4C}\int_{t}^{t+4C} \Phi(s) ds \le 4  (e^{-\omega t} + e^{-\omega (T-t)})[\Phi(0) + \Phi(T)], \qquad \text{where \, $\omega = \frac{\log 2}{4C}$}.
 \]
 \end{lemma}
 
 \begin{proof} Set first $\Psi(t) = \Phi(t) + \Phi(T-t)$. Then, for any $t \in [0, T/2]$,
 \[
 \int_{t}^{T-t} \Phi(T-s) ds =  \int_{t}^{T-t} \Phi(s) ds \le C[\Phi(t) + \Phi(T-t)] = C \Psi(t),
 \]
 and therefore
 \[
 \int_{t}^{T-t} \Psi(s)ds \le 2C \Psi(t).
 \]
 
 Since
 \[
 \int_{0}^{4C} \Psi(s)ds \le  \int_{0}^{T} \Psi(s)ds \le 2C \Psi(0),
 \]
 there exists by the Mean Value Theorem a value $\tau_1 \in [0,4C]$ such that
 \[
 \Psi(\tau_1) \le \frac12 \Psi(0).
 \]
 Then, since
  \[
 \int_{4C}^{8C} \Psi(s)ds \le  \int_{\tau_1}^{T-\tau_1} \Psi(s)ds \le 2C \Psi(\tau_1) \le C \Psi(0),
 \]
 there exists by the Mean Value Theorem a value $\tau_2 \in [4C,8C]$ such that
 \[
 \Psi(\tau_2) \le \frac14 \Psi(0).
 \]
 We can iterate this procedure to obtain a finite sequence of $\tau_n \in [4(n-1)C, 4nC]$ such that $\Psi(\tau_n) \le 2^{-n} \Psi(0)$, until $4nC \le T/2$, and
 \[
 \int_{4(n-1)C}^{4nC} \Psi(s)ds \le  \int_{\tau_n}^{T-\tau_n} \Psi(s)ds \le  \frac C{2^{n-2}} \Psi(0).
 \]
 
 Let now $t \in [0,T/2]$ and be $n$ such that $t \in [4(n-1)C, 4nC)$. If $4(n+1)C \le T/2$, then
\begin{equation}\label{estpsi}
 \int_t^{t+4C} \Psi(s) ds \le \int_{4(n-1)C}^{4nC} \Psi(s)ds + \int_{4nC}^{4(n+1)C} \Psi(s)ds \le  \frac {8C}{2^n} \Psi(0) \le 8Ce^{-\omega t}\Psi(0), \qquad \omega = \frac{\log 2}{4C},
 \end{equation}
which yields
\[
 \int_t^{t+4C} \Phi(s) ds \le 8Ce^{-\omega t} [\Phi(0) + \Phi(T)].
\]
If $4(n+1)C > T/2$ we use that $\int_t^{t+4C} \Phi(s) ds \le  \int_{\tau_{n-1}}^{T-\tau_{n-1}} \Psi(s) ds$, and conclude as before.

For $t \in [T/2,T-4C]$, we apply \eqref{estpsi} with $t \mapsto T-(t+4C)$ to get
\[
8Ce^{-\omega (T-t-4C)} [\Phi(0) + \Phi(T)] \ge  \int_{T-(t+4C)}^{T-t} \Psi(s) ds \ge \int_{T-(t+4C)}^{T-t} \Phi(T-s) ds =  \int_t^{t+4C} \Psi(s) ds,
\]
which concludes the proof.
 \end{proof}
 
 \begin{lemma}\label{poincinterp} Let $f : [t_1, t_2] \to [0,\infty)$ be Lipschitz continuous on $[t_1, t_2]$. Then,
 \[
 f^2(t) \le 2\|f'\|_{L^\infty(t_1,t_2)}(t_1-t_2) \dashint_{t_1}^{t_2}f(s) ds + \left(\dashint_{t_1}^{t_2}f(s) ds\right)^2 \qquad \forall t \in [t_1, t_2].
 \]
 \end{lemma}
 
 \begin{proof} By the Mean Value Theorem, there exists $\tau \in [t_1, t_2]$ such that $f(\tau) = \dashint_{t_1}^{t_2}f(s) ds$. Hence,
 \begin{multline*}
 f^2(t) = \int_{\tau}^t (f^2)'(s) ds +  f^2(\tau) = 2\int_{\tau}^t f(s)f'(s) ds + \left(\dashint_{t_1}^{t_2}f(s) ds\right)^2 \\ \le 2\|f'\|_{L^\infty(t_1,t_2)} \int_{t_1}^{t_2}f(s) ds + \left(\dashint_{t_1}^{t_2}f(s) ds\right)^2.
 \end{multline*}
 \end{proof}
   
We now conclude with the proof of the Poincar\'e weighted inequality stated in Theorem \ref{poincarethm}.
 
 \begin{proof}[Proof of Theorem \ref{poincarethm}] The proof is based on analogous arguments as in \cite{bakry}. 
 
 First of all we show the existence of a Lyapunov function, that is  $\phi \in C^2(\R)$, with $\phi(0)=1=\min \phi$, and $c_1|x|^2 \leq \phi(x)\leq c_2|x|^2+\tilde c$ for some $c_1, c_2,\tilde c$, which satisfies   for $r>0$,  \[-\phi''+ \phi'\bar w_x \geq \beta \phi -\gamma\chi_{B(0,r)}   \qquad \text{ in }[-\kappa^{\frac14}\pi , \kappa^{\frac14}\pi]\]    for some  constants $\beta, \gamma>0$ (depending on $r$). 
We are going to choose  $\phi=\bar w-\bar w(0)+1$. 

Using \eqref{Vkbound} and \eqref{barubarm}, we get \[  \kappa^{-\frac{1}{2}}\int_{-\pi\kappa^{\frac14}}^{\pi\kappa^{\frac14}}
   V_\kappa(y)\bar \mu(y)dy\leq \kappa^{-\frac{1}{2}}\int_{-\pi\kappa^{\frac14}}^{\pi\kappa^{\frac14}}\frac{y^2}{2} c_5e^{-c_1 y^2 } dy\leq \frac{1}{2},
   \] choosing $\kappa\geq \kappa_1$. 
   Using the fact that $V_k(x)\geq \frac{x^2}{6}$, that $\bar \lambda\leq \ell$,    and for $\kappa$ sufficiently small  \[-\bar w''+|\bar w'|^2
  = -\bar \lambda +\frac{|\bar w'|^2}{2} +V_\kappa(x) \left[1-\kappa^{-\frac{1}{2}}\int_{-\pi\kappa^{\frac14}}^{\pi\kappa^{\frac14}}
   V_\kappa(y)\bar \mu(y)dy\right]\geq -\ell +\frac{1}{2}V_\kappa(x)\geq -\ell+\frac{x^2}{12}. \] 
 Observe that  for $r<|x|<\kappa^{\frac14}\pi$, there exists $\beta=\beta(r)$ for which $\beta(\bar w(x)-\bar w(0)+1)\leq \beta (c_2x^2+c_3+c)\leq  -\ell+\frac{x^2}{12}$. 
 Now for $|x|\leq r$, it is possible to choose $\gamma=\gamma(r)>0$ such that 
  $ \gamma\geq  l-\frac{x^2}{12}+ \beta (\bar w(x)-\bar w(0)+1)$. 
  
  Now consider $f\in H^1_{\bar\mu} (-\kappa^{-\frac14}\pi, \kappa^{-\frac14}\pi)$.  Recall that $\int_{-\pi\kappa^{\frac14}}^{\pi\kappa^{\frac14} } \bar \mu(y)dy=1$ and $\bar \mu>0$. 
  First of all we observe that for all $c\in \R$ there holds 
\begin{equation}\label{proiezione} \int_{-\pi\kappa^{\frac14}}^{\pi\kappa^{\frac14} } \left(f(x)-\int_{-\pi\kappa^{\frac14}}^{\pi\kappa^{\frac14} } f(y)\mu(y) dy\right)^2 \bar \mu(x)dx
\leq \int_{-\pi\kappa^{\frac14}}^{\pi\kappa^{\frac14} } (f(x)-c)^2 \bar \mu(x)dx.  \end{equation} 
 Let us fix $r>0$, $\beta=\beta(r), \gamma=\gamma(r)$ as in the construction of the Lyapunov function. Let $c=\int_{-r}^r f(y)\bar \mu(y)dy$. 
 Then for such choice of $c$, there holds that
 \[\int_{-r}^r  (f(x)-c)^2\bar \mu(x)dx\leq C_r \int_{-r}^r f_x^2(x)\bar\mu(x)dx\] for the standard Poincar\`e inequality in the ball, with measure $\bar\mu(x)dx$. 
 Using now this inequality, the Lyapunov function,  and the fact that $\bar\mu(x)=e^{-\bar w(x)+c}$, we get
 \begin{eqnarray*} 
 && \int_{-\pi\kappa^{\frac14}}^{\pi\kappa^{\frac14} } (f(x)-c)^2 \bar \mu(x)dx\leq  \int_{-\pi\kappa^{\frac14}}^{\pi\kappa^{\frac14} }\frac{ (f(x)-c)^2 }{\beta\phi(x)} (-\phi''(x)+ \phi'(x)\bar w_x(x)+   \gamma\chi_{B(0,r)}(x) )\bar \mu(x)dx\\
 &=& \int_{-\pi\kappa^{\frac14}}^{\pi\kappa^{\frac14} }\frac{ (f(x)-c)^2 }{\beta\phi(x)} (-\phi''(x)\bar \mu(x)+ \phi'(x)\bar w_x(x)\bar\mu(x))dx +\int_{-r}^r \frac{ (f(x)-c)^2 }{\beta\phi(x)} \gamma  \bar \mu(x)dx
 \\ &\leq &   \int_{-\pi\kappa^{\frac14}}^{\pi\kappa^{\frac14} }\left(\frac{ (f(x)-c)^2 }{\beta \phi(x)}\right)_x \phi_x(x)\bar \mu(x)  dx+\frac{\gamma}{\beta} \int_{-r}^r  (f(x)-c)^2\bar \mu(x)dx\\  &\leq &  \int_{-\pi\kappa^{\frac14}}^{\pi\kappa^{\frac14} }\left[ 2\frac{ (f(x)-c)f_x(x)}{\beta \phi(x)} \phi_x(x)  - \frac{(f(x)-c)^2}{\beta \phi^2(x)}\phi_x^2(x)\right] \bar \mu(x) dx+\frac{\gamma}{\beta}C_r \int_{-r}^r f_x^2(x)\bar\mu(x)dx\\
 &\leq & \frac{1}{\beta}  \int_{-\pi\kappa^{\frac14}}^{\pi\kappa^{\frac14} }f_x^2(x) \bar\mu(x)dx+ \frac{\gamma}{\beta}C_r \int_{-r}^r f_x^2(x)\bar\mu(x)dx\leq
  \left(\frac{1+C_r\gamma}{\beta}\right) \int_{-\pi\kappa^{\frac14}}^{\pi\kappa^{\frac14} }f_x^2(x) \bar\mu(x)dx
\end{eqnarray*}   from which we conclude recalling \eqref{proiezione}. 
\end{proof}   


\medskip

\begin{flushright}
\noindent \verb"annalisa.cesaroni@unipd.it"\\
\noindent \verb"cirant@math.unipd.it"\\
Dipartimento di Matematica ``Tullio Levi-Civita" \\ Universit\`a di Padova\\
Via Trieste 63, 35121 Padova (Italy)
\end{flushright}

\end{document}